\documentstyle{article}           
\pdfoutput=1 \pdfcompresslevel=9 \pdfadjustspacing=1
\def\ifudf#1{\expandafter\ifx\csname #1\endcsname\relax}
\newif\ifpdf \ifudf{pdfoutput}\pdffalse\else\pdftrue\fi
\everydisplay={\textstyle}
\font\tbb=bbmsl10 \font\sbb=bbmsl10 scaled 700
 \font\fbb=bbmsl10 scaled 500
\newfam\bbfam \textfont\bbfam=\tbb
 \scriptfont\bbfam=\sbb \scriptscriptfont\bbfam=\fbb
\def\bb{\fam\bbfam}%
\font\tfk=eufm10 \font\sfk=eufm7 \font\ffk=eufm5
\newfam\fkfam \textfont\fkfam=\tfk
 \scriptfont\fkfam=\sfk \scriptscriptfont\fkfam=\ffk
\def\fk{\fam\fkfam}%
\def\fn[#1]{\font\TmpFnt=#1\relax\TmpFnt\ignorespaces}
\def\em{\expandafter\ifx\the\font\tensl\rm\else\sl\fi}
\def\dt{\number\day.\number\month.\number\year/\the\time}
\def\\{\hfill\break}
\def\N{{\bb N}} \def\R{{\bb R}} \def\C{{\bb C}} \def\L{{\bb L}}
\def\artanh{\mathop{\rm artanh}}
\def\Re{\mathop{\rm Re}} \def\Im{\mathop{\rm Im}}
\def\sn{\mathop{\rm sn}} \def\cn{\mathop{\rm cn}}
\def\dn{\mathop{\rm dn}} \def\am{\mathop{\rm am}}

 \def\const{{\rm const}}
\def\q{{\fk q}}
\newcount\secNo \secNo=0
\def\section #1\par{\goodbreak\vskip 3ex\noindent
 \global\advance\secNo by 1 \global\eqnNo=0
 {\fn[cmbx10 scaled 1200]#1}\vglue 1ex}
\def\href#1<#2>{\leavevmode
 \ifpdf\pdfstartlink attr {/Border [0 0 0 ]} goto name {#1}\fi
 {#2}\ifpdf\pdfendlink\fi}
\def\label@#1:#2@{\ifudf{#1}
 \expandafter\xdef\csname#1\endcsname{#2}\else
 \errmessage{label #1 already in use!}\fi}
\label@christoffel:1.1@
\label@param-ell:2.1@
\label@orthoEE:2.2@
\label@thm.EE:Lemma 2.1@
\label@orthoEM:2.3@
\label@fig.ee:Fig 1@
\label@thm.tcc:Lemma 2.2@
\label@weierform:3.1@
\label@weierdiff:3.2@
\label@elliell:3.3@
\label@intell:3.4@
\label@implell:3.5@
\label@thm.ST:Lemma 3.1@
\label@thm.CD:Lemma 3.2@
\label@fig.eed:Fig 2@
\label@thm.CE:Cor 3.3@
\label@thm.H2M:Lemma 4.1@
\label@kobayashi:4.3@
\label@hypiell:4.4@
\label@thm.CH2:Cor 4.2@
\label@fig.h2m:Fig 3@
\label@fig.h2e:Fig 4@
\label@thm.CH2E:Cor 4.3@
\label@tl-weierform:5.1@
\label@fig.h1e:Fig 5@
\label@hypihyp:5.4@
\label@hypisol:5.5@
\label@thm.CH1:Cor 5.1@
\label@fig.h1u:Fig 6@
\label@thm.CH1E:Cor 5.2@
\label@thm.laext:Lemma A.1@
\label@def.ljacobi:Def \& Cor A.2@
\label@imdg:8@
\def\@#1:#2@{\ifpdf\pdfdest name {#1} xyz\fi {#2}}
\def\:#1:{\href#1<\ifudf{#1}??\else\csname#1\endcsname\fi>}
\newcount\refNo \refNo=0
\def\refitem#1 {\global\advance\refNo by 1
 \item{\@#1:\number\refNo@.}}
\let\pleqno\eqno \newcount\eqnNo \eqnNo=0
\def\eqno#1$${\global\advance\eqnNo by 1
 \pleqno{\rm(\@#1:\number\secNo.\number\eqnNo@)}$$}
\newtoks\title \newtoks\stitle \newtoks\author
\newtoks\status \newtoks\funding
\title={Quadrics and Scherk towers}
\stitle={Quadrics \& Scherk towers}
\author={S Fujimori, U Hertrich-Jeromin, M Kokubu, M Umehara,
 K Yamada}
\status={S Fujimori et al}
\funding={This work has been partially supported by
  the Austrian Science Fund (FWF) and
  the Japan Society for the Promotion of Science (JSPS)
  through the FWF/JSPS Joint Project grant I1671-N26
  ``Transformations and Singularities''.
 Furthermore, we gratefully acknowledge partial support
  from JSPS through personal grants:
 Fujimori (B) 25800047;
 Umehara (A) 26247005;
 Yamada (C) 26400066.
 }
\ifpdf\pdfoutput=1\pdfadjustspacing=1\fi
\def\item#1{\par\leavevmode\hangindent=\parindent\hangafter=1%
 \llap{#1\enspace}\ignorespaces}
\begin{document}                  
\centerline{{\fn[cmbx10 scaled 1440]\the\title}}\vglue .2ex
\centerline{{\fn[cmr7]\the\author}}\vglue 3em plus 3ex
\centerline{\vtop{\hsize=.8\hsize{\bf Abstract.}\enspace
 We investigate the relation between quadrics and their
 Christoffel duals on the one hand, and certain zero mean
 curvature surfaces and their Gauss maps on the other hand.
 To study the relation between timelike minimal surfaces
 and the Christoffel duals of $1$-sheeted hyperboloids
 we introduce para-holomorphic elliptic functions.
 The curves of type change for real isothermic surfaces
 of mixed causal type turn out to be aligned with the real
 curvature line net.
}}\vglue 2em
\centerline{\vtop{\hsize=.8\hsize{\bf MSC 2010.}\enspace
{\it 53C42\/}, {\it 53A10\/}, 53A30, 37K25, 37K35
}}\vglue 1em
\centerline{\vtop{\hsize=.8\hsize{\bf Keywords.}\enspace
 isothermic surface; Christoffel transformation;
 minimal surface; maximal surface; saddle tower;
 Scherk surface; Karcher saddle tower;
 central quadric; ellipsoid; hyperboloid;
 timelike surface; causal type.
}}\vglue 3em

\section 1. Introduction

Though quadrics belong to the most thoroughly investigated
surfaces there are still mysterious aspects in their geometry.
In particular, quadrics can reasonably be studied in a variety
of ambient geometries
---
M\"obius or conformal geometry {\em not\/} being one of these
geometries since M\"obius transformations generally do not
map a quadric to another quadric.
Nevertheless, quadrics belong to the (M\"obius geometric)
class of isothermic surfaces,
and display rather interesting and surprising features
in this context, which are far from being understood.

This paper aims to shed some light on some of these,
seamingly incongruous, features.

Minimal surfaces in Euclidean space admit, away from umbilics,
conformal curvature line parameters:
this characterizes minimal surfaces as
{\em isothermic surfaces\/},
that is, surfaces that are ``capable of division into
infinitesimal squares by means of their curves of curvature'',
cf [\:ca72:].
Motivated by the existence of the $1$-parameter family
of associated minimal surfaces for a given minimal surface,
that are isometric and feature parallel tangent planes,
Christoffel classified those surfaces that admit a non-trivial
partner surface that is conformally related to the first by
parallel tangent planes,
cf [\:ch67:]:
apart from associated minimal surfaces these are precisely
the isothermic surfaces, which admit a partner surface that
is generically unique up to scaling so that the relation is
orientation reversing:
we will refer to this partner surface of an isothermic surface
as its {\em Christoffel dual\/}.
Given an isothermic surface $z\mapsto x(z)$ parametrized by
conformal curvature line parameters $z=u+iv$,
that is,
$$
  (x_z,x_z) = 0
   \enspace{\rm and}\enspace
  \det(x_{uv},x_u,x_v) = 0,
$$
its Christoffel dual $z\mapsto x^\ast(z)$ may be obtained
by integrating {\em Christoffel's equations\/}
$$
  x^\ast_z = {1\over(x_z,x_{\bar z})}\,x_{\bar z},
\eqno christoffel$$
cf (\:weierdiff:),
see [\:ch67:, IV] or [\:imdg:, \S5.2.1].
Any minimal surface yields an example,
with its Gauss map providing the Christoffel dual
---
the reconstruction of the minimal surface from its Gauss image
is essentially the Enneper-Weierstrass representation,
cf [\:imdg:, \S5.3.12].

It is well known that quadrics are isothermic surfaces,
see [\:ca72:],
and their Christoffel duals were determined in [\:re86:],
in terms of the usual elliptic coordinates.
A formulation using Jacobi elliptic functions,
cf [\:imdg:, \S5.2.21],
allows to study the global behaviour of the ellipsoid
as an isothermic surface,
in particular,
of its Christoffel dual and a common polarization
---
a holomorphic quadratic differential whose existence yields
another characterization of isothermic surfaces,
see [\:imdg:, \S5.2.12] and [\:sm04:, Sect 4].

The present paper was motivated by the observation that the
Christoffel dual of a tri-axial ellipsoid is the affine image
of a Scherk tower,
cf [\:imdg:, \S5.2.21 Footnote 21],
see [\:sc35:] and [\:ni89:, \S83(41)].
As the notions of isothermic surface and Christoffel duality
are M\"obius geometric and Euclidean notions, respectively,
this feature of the ellipsoid and its Christoffel dual
appear to be a highly unlikely coincidence:
the principal aim of this paper is
to shed further light on this coincidence and
to obtain a better understanding of the reasons behind it.

As a key result we derive \:thm.CD:,
which presents an approach to understand the phenomenon,
and which allows us to easily derive similar results
for $2$- and $1$-sheeted hyperboloids in Minkowski space,
cf [\:ma05:] and [\:ma04:, Prop 3.1].
As the Enneper-Weierstrass representation for minimal surfaces
in Euclidean space provides a method to explicitely determine
the Christoffel dual of a tri-axial ellipsoid,
so do Kobayashi's and Konderak's Weierstrass type
representations for maximal and timelike minimal surfaces in
Minkowski space to find the Christoffel duals of hyperboloids,
see [\:ko83:] resp [\:ko99:].
To investigate timelike minimal surfaces we derive some
results on para-holomorphic functions,
in particular,
we introduce para-complex analogues of the Jacobi elliptic
functions in \:def.ljacobi:,
cf [\:tkl14:].

As we work in Minkowski space, we obtain isothermic surfaces
that change causal type by affine transformations:
an interesting feature of these surfaces is that the lines
of separation between the space- and timelike parts of such
a surface follows the curvature line net,
see \:thm.tcc:.
Indeed, part of our investigations is independent of the
existence of and relation to a zero mean curvature surface:
we obtain explicit representations of Christoffel duals
in Euclidean as well as Minkowski ambient geometries,
see \:thm.CE:,
\:thm.CH2:, \:thm.CH2E:,
\:thm.CH1: and \:thm.CH1E:.

Though we restrict ourselves to quadrics in Minkowski space
that are aligned to the timelike axis of the ambient space,
this restriction is not essential:
the same methods will lead to similar results if a surface
is in a more general position as the nature of the occurring
differential equation will not change
---
however, computations and formulas will be less transparent.

{\it Acknowledgements.\/}
This work would not have been possible without the valuable
and enjoyable discussions with B Springborn and E Tjaden
about the subject more than a decade ago;
further we would like to thank
A Honda,
M Pember
for fruitful more recent discussions around the subject.

\the\funding

\section 2. The ellipsoid

To set the scene we discuss ellipsoids in Euclidean space,
cf [\:imdg:, \S5.2.21]:
as we wish to establish a relation with minimal surfaces,
we seek a curvature line parametrization in terms of a
meromorphic function, that is, in terms of a complex
variable.

To this end we adopt a new method to determine a suitable
curvature line parametrization of an ellipsoid,
based on two elementary observations:
{\parindent=2em
\item{$\bullet$}
 curvature line parameters $(u,v)$ of a surface
 $x:\Sigma\to\R^3$
 can be characterized as orthogonal conjugate parameters,
 and
\item{$\bullet$}
 the notion of conjugate parameters is an affine notion,
 in particular, independent of a choice of an ambient metric
 and invariant under affine transformations of $\R^3$.
\par}
If $(.,.)$ now denotes a non-degenerate inner product on $\R^3$
then $(x,x)\equiv\const$ implies that $x\perp dx$,
hence
$$
  0 = (x,x_v)_u = (x_u,x_v) + (x,x_{uv}),
$$
that is, $(u,v)$ are conjugate parameters if and only if
they are orthogonal.

In particular, conjugate parameters of the standard round
sphere $S^2\subset\R^3$ can be characterized by orthogonality
with respect to the induced metric.
Clearly, these also qualify as ``curvature line parameters'',
i.e., orthogonal conjugate parameters, on $S^2$.

In order to obtain curvature line coordinates on an ellipsoid
we parametrize the $2$-sphere conformally over a Riemann
surface $\Sigma$ and post-compose by an affine transformation,
more specifically,
$$
  \alpha x:\Sigma\to\R^3,
   \enspace{\rm where}\enspace\cases{
  x:\Sigma\to S^2 & is conformal and \cr
  \alpha:\R^3\to\R^3 & scales the axes by $a,b,c>0$. \cr}
$$
Writing $x$ in terms of a meromorphic function
$y:\Sigma\to\C\cup\{\infty\}$ we seek conditions on $y$
so that,
in terms of suitable holomorphic coordinates $z=u+iv:\Sigma\to\C$
on the Riemann surface $\Sigma$,
$$
  \alpha x
  = {1\over 1+|y|^2}\,\pmatrix{
    2a\Re y\cr 2b\Im y\cr c\,(1-|y|^2)\cr}
  = {1\over 1+y\bar y}\,\pmatrix{\hfill
    a\,(y+\bar y) \cr -ib\,(y-\bar y) \cr c\,(1-y\bar y) \cr }:
  \Sigma \to \R^3 \subset \C^3
\eqno param-ell$$
yields an orthogonal, hence curvature line parametrization.
In terms of $z=u+iv$, orthogonality of the parameter lines
$v=\const$ and $u=\const$ is expressed by the fact that the
$z$-derivative
$$
  \Im((\alpha x)',(\alpha x)') = 0,
   \enspace{\rm where}\enspace
  (\alpha x)'
  = (\alpha x)_z
  = {1\over 2}( (\alpha x)_u - i\,(\alpha x)_v )
$$
and $(.,.):\C^3\times\C^3\to\C$ denotes the bilinear extension
of the standard inner product of $\R^3$ to $\C^3$.
Excluding the case $a=b=c$, when the derivative
$$
  (\alpha x)' = {y'\over(1+y\bar y)^2}\,\pmatrix{\hfill
   a\,(1-\bar y^2) \cr -ib\,(1+\bar y^2) \cr -2c\bar y \cr}
$$
becomes isotropic,
i.e., $z\mapsto(\alpha x)(z)$ conformal,
the condition that
$z\mapsto((\alpha x)',(\alpha x)')(z)\in\R$
be real-valued reads
$$
  y'^2 = \varrho\,\{ a^2(1-y^2)^2 - b^2(1+y^2)^2 + 4c^2y^2 \},
\eqno orthoEE$$
where $z\mapsto\varrho(z)\in\R$ is a suitable real-valued
and holomorphic function, hence a (real) constant
by the Cauchy-Riemann equations $\varrho_u+i\varrho_v=0$.

Thus the function $y:\Sigma\to\C\cup\{\infty\}$ is an elliptic
function, defined on a suitable torus $\Sigma=\C/\Gamma$,
and $z$ can be considered as a globally defined coordinate
function.
For a tri-axial ellipsoid, where $a$, $b$ and $c$ are
pairwise distinct, the four branch values
$$
  y' = 0
   \enspace\Leftrightarrow\enspace
  y^2 = {1\over a^2-b^2}(\sqrt{a^2-c^2}\pm\sqrt{b^2-c^2})^2
$$
of the elliptic function $y$ yield the singularities
of the curvature line net,
that is, the four umbilics of the ellipsoid.
Note that the set of branch values is symmetric
with respect to reflections in the real and imaginary axes,
we well as with respect to inversion in the unit circle.
Hence, depending on the order of half-axis lengths, 
the branch values are all real, all purely imaginary or
all unitary.

For example, for a tri-axial ellipsoid in $\R^3$ we may
assume that $a>b>c$, without loss of generality.
We then set
$$
  p := \sqrt{a^2-b^2\over a^2-c^2}, \enspace
  q := \sqrt{b^2-c^2\over a^2-c^2} = \sqrt{1-p^2}
   \enspace{\rm and}\enspace
  r := {\sqrt{a^2-c^2}\over b};
$$
hence (\:orthoEE:) reads
$$
  y'^2
  = \varrho\,b^2r^2p^2\,
    (y^2-({1+q\over p})^2)(y^2-({1-q\over p})^2).
$$
Then the branch values $y^2=({1\pm q\over p})^2$ of $y$ are
all real and symmetric with respect to the unit circle:
this reflects the position and symmetry of the umbilics on
such a tri-axial ellipsoid, which lie on the ``equator''
ellipse in the plane orthogonal to its middle length axis.
As a constant real factor may be absorbed by a (constant)
scaling in the domain we may, without loss of generality,
assume that $1=4\varrho\,b^2r^2$ to obtain a solution
of (\:orthoEE:):
$$
  y = {1\over i}\,e^{i\am_p} = \sn_p - i\cn_p,
   \enspace{\rm where}\enspace
  \am_p:\C\to\C\cup\{\infty\}
$$
denotes the Jacobi amplitude function with module $p$,
i.e., $\am_p'=\dn_p=\sqrt{1-p^2\sn_p^2}$.
Hence a curvature line parametrization can be expressed
explicitely in terms of Jacobi elliptic functions:

\proclaim\@thm.EE:Lemma 2.1@.
Let $p,q\in(0,1)$ so that $p^2+q^2=1$ and $r\in(0,{1\over q})$;
a curvature line parametrization of a tri-axial ellipsoid is
then obtained, using the Jacobi elliptic functions
$\sn_p,\cn_p:\C\to\C\cup\{\infty\}$
with module $p$, by
$$
  u+iv=z\mapsto(\alpha x)(z) :
  = \pmatrix{\hfill
    a\sn_pu\dn_qv \cr
   -b\cn_pu\cn_qv \cr\hfill
    c\dn_pu\sn_qv \cr},
   \enspace{\sl where}\enspace\cases{
   a := \sqrt{1+r^2p^2},\cr b := 1,\cr c := \sqrt{1-r^2q^2}.\cr}
$$
Conversely, up to homothety every tri-axial ellipsoid,
with half-axes $a>b>c$,
admits such a curvature line parametrization.

Note that the same arguments apply for a Minkowski ambient
geometry:
the notion of a conjugate net is independent of metric,
hence only the orthogonality condition (\:orthoEE:) changes
by a sign
---
replacing $c$ by $ic$ in (\:orthoEE:) we obtain
$$
  y'^2 = \varrho\,\{ a^2(1-y^2)^2 - b^2(1+y^2)^2 - 4c^2y^2 \}
\eqno orthoEM$$
as the condition for $z=u+iv$ to yield a curvature line net
for the ellipsoid parametrized by (\:param-ell:).

However, there is now a distinguished (timelike) direction,
and a detailed analysis depends on how the order of
half-axes lengths of the ellipsoid relates to the
timelike direction.
In particular, the branch values of $y$ are now either all
real or all purely imaginary:
$$
  y' = 0
   \enspace\Leftrightarrow\enspace
  y^2 = {1\over a^2-b^2}(\sqrt{a^2+c^2}\pm\sqrt{b^2+c^2})^2.
$$
Also, note that an ellipsoid in Minkowski space decomposes
into three connected components,
two of which carry a positive definite induced metric
while the induced metric has Lorentz signature
on the remaining component
---
and it degenerates on the two curves separating
these components,
cf \:fig.ee:.
The umbilics necessarily lie in the spacelike part of
the ellipsoid, as is confirmed by determining the points
where the equator ellipse containing the umbilics intersects
the curve that separates the space- and timelike parts of
the ellipsoid:
when $a>b$, without loss of generality, these are given by
$$
  y^2 = {1\over a^2}\,(\sqrt{a^2+c^2}\pm c)^2;
$$
the assertion then is a consequence of the chain
of inequalities
$$
   ({\sqrt{a^2+c^2}-\sqrt{b^2+c^2}\over\sqrt{a^2-b^2}})^2
   < ({\sqrt{a^2+c^2}-c\over a})^2
   < 1
   < ({\sqrt{a^2+c^2}+c\over a})^2
   < ({\sqrt{a^2+c^2}+\sqrt{b^2+c^2}\over\sqrt{a^2-b^2}})^2.
$$

 \pdfximage width 200pt {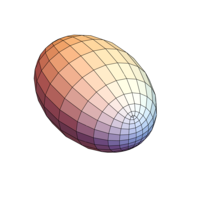}
  \xdef\elle{\the\pdflastximage}
 \pdfximage width 200pt {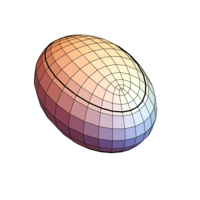}
  \xdef\ellm{\the\pdflastximage}
 \vglue -24pt\hbox to \hsize{\hfil
  \vtop{\pdfrefximage\elle}%
  \hglue -0pt
  \vtop{\pdfrefximage\ellm}%
 \hfil}\vglue -24pt\centerline{%
  {\bf\@fig.ee:Fig 1@.} A tri-axial ellipsoid with
  its curvature line net in $\R^3$ resp $\R^{2,1}$}\par

These curves of separation between the spacelike and timelike
parts of the ellipsoid follow the curvature line net,
hence may be considered as curvature lines as well:
if $(u,v)$ are orthogonal conjugate parameters
of a surface $x$, that is,
$$
  x_{uv} = \lambda x_u + \mu x_v
   \enspace{\rm and}\enspace
  (x_u,x_v) = 0,
$$
then
$$
  G_u = (x_v,x_v)_u = 2\mu\,(x_v,x_v) = 2\mu\,G
$$
satisfies a first order linear differential equations;
consequently, the function $u\mapsto G(u,v)$ vanishes
identically if it vanishes for some $u$.
Clearly, an analogous statement holds true for $E$.
Thus we have proved:

\proclaim\@thm.tcc:Lemma 2.2@.
Let $(u,v)\mapsto x(u,v)\in\R^{2,1}$ denote an
orthogonal and conjugate parametrization of a surface
patch in Minkowski space.
If the induced metric degenerates at a point $x(u,v)$
then it degenerates along a parameter curve through
this point.

\section 3. Scherk's saddle tower

We shall now see how the elliptic functions $y=\sn_p-i\cn_p$,
derived above to parametrize ellipsoids by curvature lines,
also yields (curvature line) parametrizations of Scherk's
singly periodic saddle towers,
cf [\:sc35:],
see also [\:ka89:, \S2.3.4] and [\:ni89:, \S83(41)].
Moreover, the Christoffel dual of a tri-axial ellipsoid as
an isothermic surface,
cf [\:re86:] or [\:imdg:, \S5.2.21],
can be obtained from a suitable saddle tower by an affine
transformation.

Given a meromorphic function $y:\Sigma\to\C\cup\{\infty\}$
and a polarization $\q$,
i.e., a meromorphic quadratic differential,
on the Riemann surface $\Sigma$ the Weierstrass representation
formula yields a minimal surface
$$
  x^\ast = \Re\int\pmatrix{
    1-y^2\cr i\,(1+y^2)\cr -2y\cr}\,{\q\over dy}
   \enspace{\rm with}\enspace
  x = {1\over 1+|y|^2}\,\pmatrix{
    2\Re y \cr 2\Im y \cr 1-|y|^2 \cr}
  = {1\over 1+y\bar y}\,\pmatrix{
    y+\bar y \cr {1\over i}\,(y-\bar y) \cr 1-y\bar y \cr }
\eqno weierform$$
as its Gauss map and $\q=(x^\ast_z,x_z)\,dz^2$
as its Hopf differential since
$$
  x^\ast_zdz = {\q\over(x_{\bar z},x_z)\,dz}\,x_{\bar z}.
\eqno weierdiff$$
Note that in case $\q=dz^2$, that is, when $z=u+iv$ are
conformal curvature line coordinates, the Weierstrass
formula (\:weierdiff:) simplifies to the Christoffel
formula  for the dual of an isothermic surface,
see [\:ch67:, IV] or [\:imdg:, \S5.2.1],
cf [\:imdg:, \S5.3.12]:
$$
  x^\ast_z
  = {1\over 2y'}\,\pmatrix{ 1-y^2\cr i\,(1+y^2)\cr -2y\cr}
  = {1\over(x_{\bar z},x_z)}\,x_{\bar z}
   \enspace\Leftrightarrow\enspace
  \cases{
   x^\ast_u = {2\over|x_u|^2}\,x_u, \cr
   x^\ast_v =-{2\over|x_v|^2}\,x_v. \cr}
$$

When $y:\Sigma\to\C\cup\{\infty\}$ is an elliptic function,
defined on a suitable torus $\Sigma=\C/\Gamma$,
it is straightforward to integrate the Weierstrass formula
(\:weierform:) using partial fractions.
In particular, for the function
$$
  y={1\over i}\,e^{i\am_p}=\sn_p-i\cn_p,
   \enspace{\rm satisfying}\enspace
  y'^2 = {1\over 4}\{ p^2(1+y^4)-2(1+q^2)\,y^2 \},
\eqno elliell$$
that we used in \:thm.EE: above to obtain a curvature line
parametrization of a tri-axial ellipsoid we find
$$\matrix{\hfill
  \int{1-y^2\over y'^2}\,dy
  &=&\,{2\over p}\,\artanh{2y\over p\,(1+y^2)}, \hfill\cr
  i\,\int{1+y^2\over y'^2}\,dy
  &=& {2i\over pq}\artanh{2qy\over p\,(1-y^2)} \hfill
  &=& {2\over pq}\,\arctan{2iqy\over p\,(1-y^2)}, \cr
  \int{2y\over y'^2}\,dy
  &=&\,{2\over q}\,\artanh{q\,(1+y^2)\over(1-y^2)}. \cr
}\eqno intell$$
Hence we learn (cf Appendix) that the minimal surface
$x^\ast$ of (\:weierform:) obtained from (\:elliell:)
has an implicit representation
$$
  q^2\cosh px^\ast_1
  + \cos pqx^\ast_2
  - p^2\cosh qx^\ast_3
  = 0,
\eqno implell$$
thus it yields a curvature line parametrization
of a Scherk tower,
cf [\:sc35:, \S8] and [\:ni89:, \S83(41)]:

\proclaim\@thm.ST:Lemma 3.1@.
Let $p,q\in(0,1)$ so that $p^2+q^2=1$;
a curvature line parametrization of a Scherk tower is
then obtained, using the Jacobi elliptic functions
$\sn_p,\cn_p,\dn_p:\C\to\C\cup\{\infty\}$ with module $p$, by
$$
  z\mapsto x^\ast(z) := {2\over pq}\Re\,\pmatrix{
   q\artanh{1\over p\sn_p}\cr
   \phantom{q}\arctan{q\over p\cn_p}\cr
   p\artanh{q\over i}{\sn_p\over\cn_p}\cr}(z).
$$

Thus considering $(x,x^\ast)$ as a Christoffel pair,
any affine image of the pair will form a Combescure pair:
under any affine transformation $\alpha:\R^3\to\R^3$
the common (conformal) curvature line parameters $(u,v)$
of $x$ and $x^\ast$ are turned into common conjugate parameters
of $\alpha x$ and $\alpha x^\ast$, with parallel tangents,
$$
  x^\ast_z = {1\over(x_z,x_{\bar z})}\,x_{\bar z}
   \enspace\Rightarrow\enspace
  \alpha x^\ast_z
  = {1\over(x_z,x_{\bar z})}\,\alpha x_{\bar z}.
$$
Then $\alpha x$ and $\alpha x^\ast$ are conformally related
if and only if $(u,v)$ are orthogonal coordinates,
that is, curvature line coordinates for both surfaces.
Namely, if $\alpha x^\ast_z=\varrho\,\alpha x_{\bar z}$
with a (real) function $\varrho$, then
$$
    |\alpha x^\ast_u|^2du^2
  + 2(\alpha x^\ast_u,\alpha x^\ast_v)\,dudv
  + |\alpha x^\ast_u|^2du^2
  = \varrho^2\{
    |\alpha x_u|^2du^2
  - 2(\alpha x_u,\alpha x_v)\,dudv
  + |\alpha x_u|^2du^2
    \}.
$$
In fact, we obtain a stronger result, as integrability of
$dx^\ast$ implies that $z=u+iv$ yield conjugate parameters,
since
$$
  0 = x^\ast_{uv} - x^\ast_{vu}
  = (\varrho x_u)_v + (\varrho x_v)_u
  = 2\varrho x_{uv} + \varrho_vx_u + \varrho_ux_v.
$$
Thus, using Christoffel's characterization [\:ch67:]
of a Christoffel pair of isothermic surfaces as an
orientation reversing conformal Combescure pair,
cf [\:ma05:, Main Thm] or [\:ma04:, Prop 3.1],
we obtain the following

\proclaim\@thm.CD:Lemma 3.2@.
Suppose that the differentials of two surfaces $x$ and $x^\ast$
in $\R^3$ are related by
$$
  x^\ast_z = \varrho\,x_{\bar z}
  \enspace\Leftrightarrow\enspace\cases{
   x^\ast_u = \varrho\,x_u, \cr
   x^\ast_v = -\varrho\,x_v. \cr}
$$
Then
$(u,v)$ are common conjugate parameters for $x$ and $x^\ast$.
Moreover, if $(u,v)$ are orthogonal coordinates
for some ambient metric $(.,.)$,
$$
  x_u\perp x_v
   \enspace\Leftrightarrow\enspace
  z\mapsto(x_z,x_z)(z)\in\R,
$$
then $(x,x^\ast)$ is a Christoffel pair with
common curvature line parameters $(u,v)$.

Note that, starting as above with common conformal curvature
line coordinates $(u,v)$ of a Christoffel pair $(x,x^\ast)$
to obtain a new Christoffel pair $(\alpha x,\alpha x^\ast)$
by an affine transformation via the orthogonality condition
of \:thm.CD:, the coordinates will in general not be conformal
for either surface, $\alpha x$ or $\alpha x^\ast$:
in general $(\alpha x_z,\alpha x_z)$ does not vanish,
even if it is real.
However, as $x$ satisfies a Laplace equation
$$
  x_{uv} = \lambda x_u + \mu x_v
$$
when $(u,v)$ are conjugate parameters,
we infer that
$$
  {1\over 2}(\ln{(\alpha x_u,\alpha x_u)\over(x_u,x_u)})_v
  = \mu\,\{
    {(\alpha x_u,\alpha x_v)\over(\alpha x_u,\alpha x_u)}
    - {(x_u,x_v)\over(x_u,x_u)} \} = 0
   \enspace{\rm and}\enspace
  {1\over 2}(\ln{(\alpha x_v,\alpha x_v)\over(x_v,x_v)})_u = 0
$$
as soon as $(u,v)$ are curvature line coordinates for both,
$x$ and $\alpha x$.
Consequently, conformal curvature line coordinates
$(\tilde u,\tilde v)$ for $\alpha x$ are obtained from
those of $x$ by integrating
$$
  d\tilde u = {|\alpha x_u|\over|x_u|}\,du
   \enspace{\rm and}\enspace
  d\tilde v = {|\alpha x_v|\over|x_v|}\,dv.
$$
Identification of conformal curvature line parameters for
$\alpha x$ is more involved if a reference conformal factor
is missing, that is, if $(u,v)$ were not conformal curvature
line parameters for $x$.

Combining the results from \:thm.EE:, \:thm.ST: and \:thm.CD:
we learn how a solution $y$ of (\:orthoEE:) gives rise
to a curvature line parametrization of a tri-axial ellipsoid
 on the one hand, and
to a curvature line parametrization of a Scherk tower
 on the other hand
---  hence, by \:thm.CD:,
to a curvature line parametrization of the Christoffel dual
of the tri-axial ellipsoid,
cf [\:re86:,\S6]:

 \pdfximage width 200pt {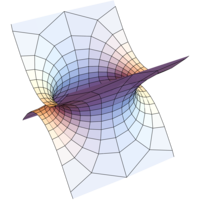}
  \xdef\elled{\the\pdflastximage}
 \pdfximage width 200pt {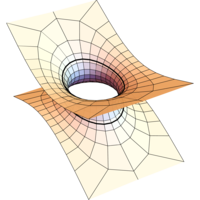}
  \xdef\ellmd{\the\pdflastximage}
 \vglue -0pt\hbox to \hsize{\hfil
  \vtop{\pdfrefximage\elled}%
  \hglue -0pt
  \vtop{\pdfrefximage\ellmd}%
 \hfil}\vglue 8pt\centerline{%
  {\bf\@fig.eed:Fig 2@.} The Christoffel duals of the
  ellipsoids in $\R^3$ resp $\R^{2,1}$ of \:fig.ee:}\par

\proclaim\@thm.CE:Cor 3.3@.
Let $p,q\in(0,1)$ so that $p^2+q^2=1$ and $r\in(0,{1\over q})$;
the Christoffel dual of the tri-axial ellipsoid with half axes
$$
  a := \sqrt{1+r^2p^2}, \enspace
  b := 1 \enspace{\sl and}\enspace
  c := \sqrt{1-r^2q^2}
$$
is an affine transformation of a Scherk tower, thus in terms
of Jacobi elliptic functions given by
$$
  z\mapsto\alpha x^\ast(z) := \Re\,\pmatrix{
   {2a\over p}\artanh{1\over p\sn_p}\cr
   {2b\over pq}\arctan{q\over p\cn_p}\cr
   {2c\over q}\artanh{q\over i}{\sn_p\over\cn_p}\cr}(z).
$$

\:thm.CD: holds also for spacelike or timelike real isothermic
surfaces in a Minkowski ambient geometry:
apart from the characterization of curvature line coordinates
as orthogonal conjugate parameters the proof of \:thm.CD:
only relies on Christoffel's formula
$$
  x^\ast_z = {1\over(x_z,x_{\bar z})}\,x_{\bar z},
$$
which holds in both the spacelike and timelike real
isothermic cases,
cf [\:ma05:, Main Thm].

However, the standard $2$-sphere $S^2\subset\R^3$ is not
totally umbilic in $\R^{2,1}$, hence inverse stereographic
projection of a holomorphic function $y$ does in general
not yield an orthogonal net on $S^2$.
In particular, the solution $y$ of (\:orthoEM:) that provides
a curvature line net of a tri-axial ellipsoid by inverse
stereographic projection followed by a suitable affine
transformation does in general not yield a curvature line
net on the $2$-sphere, before the affine transformation:
setting $a=b$ in (\:orthoEM:) confirms that a curvature
line net of $S^2$ as a non-umbilic surface of revolution
in $\R^{2,1}$ is necessarily given by its meridians and
circles of latitude.

Note that the catenoid is still the Christoffel dual of
the (Euclidean) sphere as a surface of revolution
in Minkowski space,
since the meridians and lines of latitude are orthogonal:
the Christoffel dual $(u,v)\mapsto x^\ast(u,v)$ of a surface
of revolution $(u,v)\mapsto x(u,v)$
can explicitely be computed as
$$
  x^\ast(u,v) = \pmatrix{
   -{1\over r(u)}\cos v\cr
   -{1\over r(u)}\sin v\cr
   \int{h'(u)du\over r^2(u)}\cr},
   \enspace{\rm where}\enspace
  x(u,v) = \pmatrix{ r(u)\cos v\cr r(u)\sin v\cr h(u)\cr},
$$
since
$$
  x^\ast_z(u,v) = {1\over r^2(u)}\,x_{\bar z}(u,v)
   \enspace{\rm and}\enspace
  (x_z(u,v),x_z(u,v)) = {(r'^2\pm h'^2)-r^2\over 4}(u) \in\R
$$
for the standard metrics $(.,.)$ of signature $(+,+,\pm)$
on $\R^3$.
Hence \:thm.CD: shows that the Minkowski Christoffel dual of
a surface of revolution with timelike axis coincides with its
Euclidean dual.
In particular, for $r(u)={1\over\cosh u}$ and $h(u)=\tanh u$,
the conformal curvature line parametrizations of the catenoid
and its Gauss map are obtained.

\section 4. The spacelike hyperboloid

In order to obtain a similar relation between
the Christoffel duality of a quadric and
an affine transformation of a zero mean curvature surface
in Minkowski space,
as we did for the Christoffel dual of an ellipsoid
in Euclidean space in the previous section,
we will need to focus on a totally umbilic quadric
in $\R^{2,1}$.
As we shall see, the same methodology then yields
similar results.

Thus consider inverse stereographic projection onto the
standard $2$-sheeted hyperboloid in $\R^{2,1}$ as a conformal
parametrization of the hyperboloid:
$$
  (\C\cup\{\infty\})\setminus S^1\ni y \mapsto x :
  = {1\over 1-|y|^2}\,\pmatrix{
    2\Re y\cr 2\Im y\cr 1+|y|^2\cr}
  = {1\over 1-y\bar y}\,\pmatrix{\hfill
    (y+\bar y) \cr -i\,(y-\bar y) \cr 1+y\bar y \cr }
  \in \R^{2,1} \subset \C^3.
$$
A parametrization of a general $2$-sheeted hyperboloid,
with timelike principal axis,
is then obtained as as an affine transformation of the standard
hyperboloid,
$$
  \alpha x
  = {1\over 1-|y|^2}\,\pmatrix{
    2a\Re y\cr 2b\Im y\cr c\,(1+|y|^2)\cr}
  = {1\over 1-y\bar y}\,\pmatrix{\hfill
    a\,(y+\bar y) \cr -ib\,(y-\bar y) \cr c\,(1+y\bar y) \cr }:
  \Sigma \to \R^{2,1} \subset \C^3.
\eqno param-2hyp$$

As before, the orthogonality condition
 $z\mapsto((\alpha x)',(\alpha x)')(z)\in\R$
for the conjugate parameters $(u,v)$, written as a complex
parameter $z=u+iv$, yields an elliptic differential equation
$$
  y'^2 = \varrho\,\{ a^2(1+y^2)^2 - b^2(1-y^2)^2 - 4c^2y^2 \},
\eqno orthoH2M$$
with its branch values
$$
  y^2 = -{1\over a^2-b^2}(\sqrt{a^2-c^2}\pm\sqrt{b^2-c^2})^2
$$
determining the umbilics of the $2$-sheeted hyperboloid in
Minkowski space $\R^{2,1}$.
Note that (\:orthoH2M:) may be obtained from the equation
(\:orthoEE:) for a Euclidean ellipsoid
by replacing $y$ by $iy$ and
a sign change of $\varrho$.
However,
in contrast to the case of an ellipsoid in Euclidean space,
the hyperboloid has a distinguished (timelike) axis and
various geometric configurations may occur,
as an elementary case study reveals.

{\parindent=2em
\item{(i)} $a,b\geq c$:
 in this case, the $2$-sheeted hyperboloid is spacelike
 and, as long as the half-axis lengths are pairwise distinct,
 the branch values of $y$ are either real or imaginary and
 yield the four umbilics of the $2$-sheeted hyperboloid,
 with their symmetries reflected by the symmetries of the
 branch values of $y$, as in the case of the ellipsoid,
 cf Fig 5.
 If $a=b$ the hyperboloid is a surface of revolution with only
 two umbilics, given by branch values at $y=0$ and $y=\infty$.
 If, on the other hand, $a=c$ or $b=c$ then the branch values
 of $y$ become unitary, reflecting the absence of umbilics
 in this case.

\item{(ii)} $a,b<c$:
 now the hyperboloid has mixed causal type,
 each component consists of a (compact) spacelike part and
 an (annular) timelike part.
 Unless $a=b$ and we obtain a surface of revolution,
 the hyperboloid has four umbilics in its spacelike part,
 given by the real or imaginary branch values of $y$,
 cf Fig 6.
 Note that, by \:thm.tcc:, the lines where the causal type
 of the hyperboloid changes are aligned with
 the curvature line net.

\item{(iii)} $a<c\leq b$ or $a\geq c>b$:
 in these cases, the hyperboloid has also mixed causal type,
 but now both components decompose into one spacelike and
 two timelike parts,
 cf Fig 7,
 and the branch values of $y$ become unitary,
 reflecting the fact that these hyperboloids have no umbilics.
\par}

Thus, in analogy to \:thm.EE: for curvature line parameters
of Euclidean ellipsoids, we obtain curvature line parameters
for $2$-sheeted hyperboloids in Minkowski space in terms of
Jacobi elliptic functions,
based on the same methods as in the ellipsoid case:

\proclaim\@thm.H2M:Lemma 4.1@.
Any tri-axial $2$-sheeted hyperboloid in Minkowski space
admits, up to homothety, a curvature line paramtrization
using the Jacobi elliptic functions
$\cn_p,\sn_p:\C\to\C\cup\{\infty\}$,
$$
  u+iv=z\mapsto\alpha x(z) : 
  = {1\over\dn_pu\sn_qv}\,\pmatrix{
    a\cn_pu\cn_qv \cr b\sn_pu\dn_qv \cr c \cr},
   \enspace{\sl where}\enspace\cases{
  a := \sqrt{1+r^2p^2},\cr b := 1,\cr c := \sqrt{1-r^2q^2}\cr}
$$
and $p,q,r\in\R\cup i\R$ with $p^2+q^2=1$ satisfy
one of the following:
$$
  {\rm(i)}\enspace
   p,q\in(0,1), r\in(0,{1\over q}); \quad
  {\rm(ii)}\enspace
   p\in i\R, q>1, {r\over i}\in(0,{1\over q}); \quad
  {\rm(iii)}\enspace
   p>1, q\in i\R, {r\over i}\in(0,{1\over p}).
$$

As \:thm.CD:, which served as the key step in characterizing
the Christoffel dual of a Euclidean ellipsoid as an affine
image of a Scherk saddle tower,
did not depend on the signature of the ambient metric
an analogous result is obtained for the Christoffel duals
of $2$-sheeted hyperboloids in Minkowski space,
cf \:thm.CE:.

Employing the Weierstrass representation of [\:ko83:, Thm 1.1]
for maximal surfaces in Minkowski space,
$$
  x^\ast = \Re\int\pmatrix{
    1+y^2\cr i\,(1-y^2)\cr 2y\cr}\,{\q\over dy}
   \enspace{\rm with}\enspace
  x
  = {1\over 1-y\bar y}\,\pmatrix{
    y+\bar y \cr {1\over i}\,(y-\bar y) \cr 1+y\bar y \cr }
\eqno kobayashi$$
as the Gauss map and $\q=dz^2$ as a polarization of
the (universal cover of the) underlying Riemann surface
$\Sigma=\C$ yields again Christoffel's formula for the dual
of a (spacelike) isothermic surface,
$$
  x^\ast_z
  = {1\over 2y'}\,\pmatrix{ 1+y^2\cr i\,(1-y^2)\cr 2y\cr}
  = {1\over(x_{\bar z},x_z)}\,x_{\bar z}
   \enspace\Leftrightarrow\enspace
  \cases{
   x^\ast_u = {2\over|x_u|^2}\,x_u, \cr
   x^\ast_v =-{2\over|x_v|^2}\,x_v. \cr}
$$

Thus replacing $y$ by $iy$ in (\:elliell:), to obtain a
solution of (\:orthoH2M:) with $\varrho=-{1\over 4(a^2-c^2)}$
from that of (\:orthoEE:), we arrive at
$$
  y=e^{i\am_p}=\cn_p+i\sn_p,
   \enspace{\rm satisfying}\enspace
  y'^2 = -{1\over 4}\{ p^2(1+y^4)+2(1+q^2)\,y^2 \}.
\eqno hypiell$$
Using (\:intell:) we hence integrate (\:kobayashi:),
$$\matrix{\hfill
  \int{1+y^2\over y'^2}\,dy
  &=&\,-{2\over p}\,\arctan{2y\over p\,(1-y^2)}, \hfill\cr
  i\,\int{1-y^2\over y'^2}\,dy
  &=& -{2i\over pq}\arctan{2qy\over p\,(1+y^2)} \hfill
  &=&  {2\over pq}\,\artanh{2qy\over ip\,(1+y^2)}, \cr
  \int{2y\over y'^2}\,dy
  &=&\,\phantom{-}
       {2\over q}\,\artanh{q\,(1-y^2)\over(1+y^2)}. \cr
}\eqno int2hyp$$
As the timelike $x_3$-axis is geometrically distinguished
we obtain (cf Appendix) two qualitatively different implicit
representations of the maximal surfaces $x^\ast$ obtained from
(\:hypiell:),
cf [\:fkk15:, Sect 3]:
$$\matrix{
  q^2\cos px^\ast_1
  + \cosh pqx^\ast_2
  - p^2\cosh qx^\ast_3
  = 0
  &{\rm if}& p<1 \enspace{\rm and}\enspace q\in(0,1); \hfill\cr
  q^2\cos px^\ast_1
  +\kern 6pt \cos{pqx^\ast_2\over i}
  -\kern 6pt p^2\cos{qx^\ast_3\over i}
  = 0
  &{\rm if}& p>1 \enspace{\rm and}\enspace q\in i\R. \hfill\cr
}\eqno impl2hyp$$
For the tri-axial $2$-sheeted hyperboloid this yields another
``permutability theorem'', intertwining Christoffel duality
and affine transformation:

\proclaim\@thm.CH2:Cor 4.2@.
The Christoffel dual of a tri-axial $2$-sheeted hyperboloid
in Minkowski space with half axes lengths
$$
  {\rm(i)}\enspace a>b>c, \quad
  {\rm(ii)}\enspace c>a>b \quad{\sl or}\quad
  {\rm(iii)}\enspace b>c>a \quad
$$
is the affine image of a maximal surface,
more precisely,
up to homothety it is given by
$$
  z\mapsto\alpha x^\ast(z) := \Re\,\pmatrix{
   {2a\over p}\arctan{1\over ip\sn_p}\cr
   {2b\over pq}\artanh{q\over ip\cn_p}\cr
   {2c\over q}\artanh{q\over i}{\sn_p\over\cn_p}\cr}(z)
  \enspace{\sl with}\enspace\cases{
   p := \sqrt{a^2-b^2\over a^2-c^2}, \cr
   q := \sqrt{b^2-c^2\over a^2-c^2}. \cr}
$$

 \pdfximage width 140pt {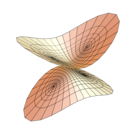}
  \xdef\hypmA{\the\pdflastximage}
 \pdfximage width 140pt {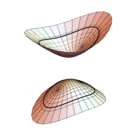}
  \xdef\hypmB{\the\pdflastximage}
 \pdfximage width 100pt {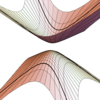}
  \xdef\hypmC{\the\pdflastximage}
 \pdfximage width 140pt {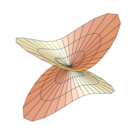}
  \xdef\hypmAd{\the\pdflastximage}
 \pdfximage width 140pt {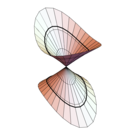}
  \xdef\hypmBd{\the\pdflastximage}
 \pdfximage width 140pt {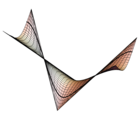}
  \xdef\hypmCd{\the\pdflastximage}
 \vglue -0pt\hbox to \hsize{\hfil
  \vtop{\pdfrefximage\hypmA}%
  \hglue -0pt
  \vtop{\pdfrefximage\hypmB}%
  \hglue -0pt
  \hglue 20pt\vtop{\vglue 16pt\pdfrefximage\hypmC}\hglue 20pt%
 \hfil}\vglue -8pt\hbox to \hsize{\hfil
  \vtop{\pdfrefximage\hypmAd}%
  \hglue -0pt
  \vtop{\pdfrefximage\hypmBd}%
  \hglue -0pt
  \vtop{\pdfrefximage\hypmCd}%
 \hfil}\vglue 8pt\centerline{%
  {\bf\@fig.h2m:Fig 3@.} Three types of $2$-sheeted
   hyperboloids and their Christoffel duals in $\R^{2,1}$}\par

Note that the hyperboloid as well as its dual are spacelike
in case (i), whereas both surfaces change causal type
in the cases (ii) and (iii),
even though the employed maximal surfaces are spacelike,
cf \:fig.h2m:.
On their timelike parts the surfaces are real isothermic
in the sense of [\:ma05:, Sect 2.1],
that is, admit Lorentz conformal curvature line parameters,
and the curvature line nets on the spacelike and timelike parts
of the surfaces extend across the lines of causal type change,
which follow the curvature line net by \:thm.tcc:.
In case (iii) the Christoffel dual surface has a translational
period, cf \:fig.h2m:,
and the surface can be extended to a triply periodic surface
in Minkowski space,
as is also seen from the implicit representation (\:impl2hyp:).

As the $2$-sheeted hyperboloid is not totally umbilic
in Euclidean space,
the construction of the Christoffel dual as an affine image
of a minimal surface fails for a hyperboloid in a Euclidean
ambient geometry.
However, the only failure of the construction turns out
to be the minimality:
integrating
$$
  x^\ast_z
  = {1\over 2y'}\,\pmatrix{ 1+y^2\cr i\,(1-y^2)\cr 2y\cr}
  = {2|y'|^2\over(1-|y|^2)^2}\,x_{\bar z}
   \enspace{\rm with}\enspace
  x = {1\over 1-|y|^2}\,\pmatrix{
    2\Re y \cr 2\Im y \cr 1+|y|^2 \cr}
$$
we obtain $x^\ast$ as the (harmonic) real part of a holomorphic
curve in $\C^3$ that is not a null curve though,
hence does not yield a minimal surface $x^\ast$ in $\R^3$.
As the pair $(x,x^\ast)$ satisfies the assumption of \:thm.CD:
we just wheel out the orthogonality condition
 $z\mapsto((\alpha x)',(\alpha x)')(z)\in\R$
for some affine image $\alpha x$ of $x$ to obtain
a Christoffel pair $(\alpha x,\alpha x^\ast)$,
with common curvature line coordinates given by a solution $y$
of the elliptic differential equation
$$
  y'^2 = \varrho\,\{ a^2(1+y^2)^2 - b^2(1-y^2)^2 + 4c^2y^2 \},
\eqno orthoH2E$$
with real or imaginary branch values given by
$$
  y^2 = -{1\over a^2-b^2}(\sqrt{a^2+c^2}\pm\sqrt{b^2+c^2})^2.
$$
As long as $a\neq b$,
that is, the hyperboloid is not a surface of revolution,
the four branch values are symmetric with respect to
the origin as well as to the unit circle,
reflecting the symmetry of the four umbilics of
the $2$-sheeted hyperboloid, cf \:fig.h2e:.

 \pdfximage width 200pt {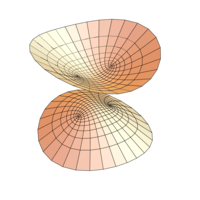}
  \xdef\hypeA{\the\pdflastximage}
 \pdfximage width 200pt {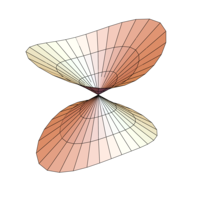}
  \xdef\hypeAd{\the\pdflastximage}
 \vglue -0pt\hbox to \hsize{\hfil
  \vtop{\pdfrefximage\hypeA}%
  \hglue -0pt
  \vtop{\pdfrefximage\hypeAd}%
 \hfil}\vglue -8pt\centerline{%
  {\bf\@fig.h2e:Fig 4@.} The Euclidean $2$-sheeted hyperboloid
  and its Christoffel dual, cf \:fig.h2m:}\par

More precisely, assuming without loss of generality $a>b$
for a non-rotational $2$-sheeted hyperboloid,
and with $0=1+4\varrho(a^2+c^2)$,
we obtain
$$
  y = \cn_p+i\sn_p = e^{i\am_p}
   \enspace{\rm with}\enspace
  p = \sqrt{a^2-b^2\over a^2+c^2} \in (0,1)
$$
as a solution of (\:orthoH2E:),
completely analogous to (\:orthoH2M:).
Consequently, a curvature line parametrization of the
Christoffel dual of a $2$-sheeted hyperboloid in Euclidean
space is obtained by precisely the same formulas as in
Minkowski space,
in \:thm.CH2:,
cf [\:re86:,\S8]:

\proclaim\@thm.CH2E:Cor 4.3@.
The Christoffel dual of a non-rotational $2$-sheeted
hyperboloid in Euclidean space with half axes lengths
$a,b,c>0$, where
$
  1+({x_1\over a})^2+({x_2\over b})^2=({x_3\over c})^2
$
and $a>b$ without loss of generality, admits a curvature line
parametrization $z\mapsto\alpha x^\ast(z)$ by Jacobi elliptic
functions:
$$
  z\mapsto\alpha x^\ast(z) := \Re\,\pmatrix{
   {2a\over p}\arctan{1\over ip\sn_p}\cr
   {2b\over pq}\artanh{q\over ip\cn_p}\cr
   {2c\over q}\artanh{q\over i}{\sn_p\over\cn_p}\cr}(z)
  \enspace{\sl with}\enspace\cases{
   p := \sqrt{a^2-b^2\over a^2+c^2}, \cr
   q := \sqrt{b^2+c^2\over a^2+c^2}. \cr}
$$

\section 5. The timelike hyperboloid

In the previous section we investigated the dual of a tri-axial
$2$-sheeted hyperboloid as an affine transform of a maximal
surface that was obtained as a dual of the standard totally
umbilic $2$-sheeted hyperboloid in Minkowski space $\R^{2,1}$.
However, the Lorentz sphere yields another totally umbilic
quadric, thus provides a starting point for a similar
investigation of the Christoffel transforms
of $1$-sheeted hyperboloids.

As the $1$-sheeted hyperboloid
  $S^{1,1}=\{x\in\R^{2,1}\,|\,(x,x)=1\}$
is a timelike surface in Minkowski space Christoffel duality
leads to Konderak's Weierstrass-type representation
[\:ko99:, Thm 3.3]
for minimal timelike surfaces in Minkowski space:
instead of holomorphic functions this representation
employs para-holomorphic functions,
defined on the algebra of Lorentz-numbers,
$$
  \L = \{y_1+jy_2\,|\,y_1,y_2\in\R\},
   \enspace{\rm where}\enspace
  j^2 = 1,
$$
cf [\:co49:] or [\:haha04:].
Note that $\L$ can be considered as the Clifford algebra
of the real line with its usual operations.
More specifically, we consider the inverse stereographic
projection
$$
  \L\setminus S^1\ni y \mapsto x :
  = {1\over 1-\|y\|^2}\,\pmatrix{
    2\Re y\cr 2\Im y\cr 1+\|y\|^2\cr}
  = {1\over 1-y\bar y}\,\pmatrix{
    y+\bar y \cr {1\over j}\,(y-\bar y) \cr 1+y\bar y \cr }
  \in S^{1,1}\subset\R^{2,1}
$$
into the Minkowski space of signature $(-,+,+)$,
where $\bar y$ and $\|y\|$ denote the conjugate
resp the modulus of a Lorentz number $y\in\L$,
$$
  \|y\|^2 = y\bar y = (y_1+jy_2)(y_1-jy_2) = y_1^2-y_2^2
   \enspace{\rm and}\enspace
  \cases{ \Re(y_1+jy_2) = y_1, \cr \Im(y_1+jy_2) = y_2, \cr}
$$
as usual.
With a para-holomorphic function $y:\L\to\L$, $z\mapsto y(z)$,
an analogue of the Enneper-Weierstrass representation formula
then yields a conformal curvature line parametrization
$$
  x^\ast = -\Re\int\pmatrix{
    1+y^2\cr j\,(1-y^2)\cr 2y\cr}\,{dz\over y'}
\eqno tl-weierform$$
of a timelike minimal surface with Gauss map $x$,
cf [\:ko99:, Thm 3.3]:
as the Christoffel formula
$$
  x^\ast_z
  = -{1\over 2y'}\,\pmatrix{ 1+y^2\cr j\,(1-y^2)\cr 2y\cr}
  = {1\over(x_{\bar z},x_z)}\,x_{\bar z}
   \enspace\Leftrightarrow\enspace
  \cases{
   x^\ast_u = {2\over|x_u|^2}\,x_u, \cr
   x^\ast_v = {2\over|x_v|^2}\,x_v \cr}
$$
holds with $z=u+jv$,
by the para-complex version of the Cauchy-Riemann equations
$$
  0 = 2y_{\bar z}
  = ((y_1)_u-(y_2)_v) + j\,((y_2)_u-(y_1)_v),
$$
we deduce that $x\perp dx^\ast$ and,
with \:thm.CD:,
that $(u,v)$ are conjugate parameters for $x^\ast$
which are orthogonal as they are conformal,
$$
  (x_z,x_z)
  = {1\over 4}\{ (|x_u|^2+|x_v|^2) + 2j\,(x_u,x_v) \}
  = 0.
$$
Note that $(x,x^\ast)$ forms a pair of {\em real isothermic\/}
surfaces, in the sense of [\:ma05:, Sect 2]
as the surfaces have real curvature directions:
in general, the shape operator of a timelike surface
in Minkowski space may have complex conjugate eigendirections,
or not diagonalize at all.

Wheeling out the orthogonality condition for the coordinates
of an affine transform
$$
  \alpha x
  = {1\over 1-\|y\|^2}\,\pmatrix{
    2a\Re y\cr 2b\Im y\cr c\,(1+\|y\|^2)\cr}
  = {1\over 1-y\bar y}\,\pmatrix{\hfill
    a\,(y+\bar y) \cr jb\,(y-\bar y) \cr c\,(1+y\bar y) \cr }:
  \Sigma \to \R^{2,1} \subset \L^3,
\eqno param-1hyp$$
as in the spacelike cases,
we arrive now at the Lorentzian ordinary differential equation
$$
  y'^2 = -\varrho\,\{ a^2(1+y^2)^2 - b^2(1-y^2)^2 - 4c^2y^2 \}
\eqno orthoH1M$$
with a real function $\varrho$ as the condition for $z=u+jv$
to yield a (real) curvature line net for the hyperboloid
parametrized by (\:param-1hyp:),
cf (\:orthoEE:) and (\:orthoH2M:).
Then, similar to the complex case, the para-complex
Cauchy-Riemann equations imply that the real function
$\varrho$ is, in fact, a real constant,
since $y_{\bar z}=0$ imples $\varrho_{\bar z}=0$.

 \pdfximage width 200pt {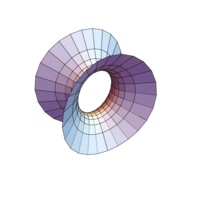}
  \xdef\hypeD{\the\pdflastximage}
 \pdfximage width 200pt {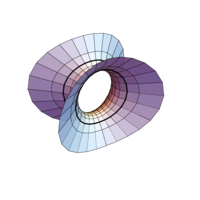}
  \xdef\hypmD{\the\pdflastximage}
 \vglue -30pt\hbox to \hsize{\hfil
  \vtop{\pdfrefximage\hypmD}%
  \hglue -0pt
  \vtop{\pdfrefximage\hypeD}%
 \hfil}\vglue -28pt\centerline{%
  {\bf\@fig.h1e:Fig 5@.} Parametrizations of a $1$-sheeted
  hyperboloid in $\R^3$ resp $\R^{2,1}$
  obtained from (\:hypihyp:)}\par

Note though that the apparent similarity with the cases
previously examined is deceptive:
in contrast to the elliptic equation (\:orthoH2M:),
or (\:orthoEE:),
we are now dealing with a para-holomorphic function $y$
and its para-complex derivative in (\:orthoH1M:).
However, as the differential equations are real,
we may employ a real solution of (\:orthoH2M:) to obtain
a real solution of (\:orthoH1M:) that we may then extend
to a para-holomorphic solution.
To make this idea more tangible consider our above solution
$$
  y(z) = e^{i\am_pz} = \cn_pz + i\sn_pz,
   \enspace{\rm where}\enspace
  p = \sqrt{a^2-b^2\over a^2-c^2},
$$
of the elliptic equation (\:orthoH2M:),
$$
  y'^2 = -{1\over 4}\{ p^2(y^4+1) + 2(1+q^2)\,y^2 \}
  = -{1\over 4(a^2-c^2)}
    \{ a^2(1+y^2)^2 - b^2(1-y^2)^2 - 4c^2y^2 \}.
$$
We may then replace $y(z)$ by $y(iz)$ to obtain a real solution
of (\:orthoH1M:) with a suitable choice of the constant,
$\varrho=-{1\over 4(a^2-c^2)}$,
by restriction of the complex solution to real variables:
$$
  y(z) = e^{i\am_piz}
  = {1-\sn_qz\over\cn_qz}
  = {\cn_qz\over 1+\sn_qz}
   \enspace{\rm solves}\enspace
  y'^2 = {1\over 4}\{ p^2(1+y^4) + 2(1+q^2)\,y^2 \}.
\eqno hypihyp$$
Note that the conversion formulae for reciprocal and imaginary
moduli do not affect reality of the Jacobi elliptic functions,
hence this approach does not depend on an order of the half
axes lengths.

We then extend this real analytic solution
(uniquely, cf \:thm.laext:)
to a Lorentz-analytic solution
$$
  y(u+jv) = y(u+v)\,{1+j\over 2} + y(u-v)\,{1-j\over 2}
  = {\cn_qu-j\sn_qv\dn_qu\over\cn_qv+\sn_qu\dn_qv}.
\eqno hypisol$$
Note that the same arguments that prove the corresponding
properties of the $\L$-Jacobi functions in \:def.ljacobi:
also show that this $\L$-analytic extension indeed solves
the para-holomorphic differential equation (\:orthoH1M:).
The Lorentz-analytic function $y$ then yields a curvature
line parametrization of (part of) the $1$-sheeted hyperboloid,
cf \:fig.h1e:.

As the differential equation (\:hypihyp:) that $y$ now
satisfies only differs by a sign from (\:hypiell:)
we read the integrals in (\:tl-weierform:) off
(\:int2hyp:),
$$\matrix{\hfill
  \int{1+y^2\over y'^2}\,dy
  &=&\phantom{-}\,{2\over p}\,
      \arctan{2y\over p\,(1-y^2)}, \hfill\cr
  j\,\int{1-y^2\over y'^2}\,dy
  &=&\phantom{-}
      {2j\over pq}\arctan{2qy\over p\,(1+y^2)}, \hfill\cr
  \int{2y\over y'^2}\,dy
  &=&\,-{2\over q}\,\artanh{q\,(1-y^2)\over(1+y^2)}. \cr
}\eqno int1hyp$$
In case of a complex elliptic modulus, $p\in i\R$,
the terms may be turned into purely para-complex form
by replacing $\arctan$ by $\artanh$,
cf Appendix.
Again we obtain two qualitatively different implicit
representations for the timelike minimal surfaces
$x^\ast$ given by (\:hypihyp:),
where the $x_1$-axis is now the distinguished axis:
$$\matrix{\hfill
  q^2\cos\Re pz_1
  +\kern 5pt \cos\Re pqz_2
  - p^2\cosh\Re qz_3
  = 0
  &{\rm if}& p\in(0,1) \enspace{\rm and}\enspace q\in(0,1); \hfill\cr
  q^2\cosh\Re{pz_1\over i}
  + \cosh\Re{pqz_2\over i}
  + p^2\cosh\Re qz_3
  = 0
  &{\rm if}& p\in i\R \enspace{\rm and}\enspace q\in(1,\infty).
   \hfill\cr
}\eqno impl1hyp$$

Thus, analogous to \:thm.CH2:, we obtain a ``permutability''
result for the Christoffel dual of a $1$-sheeted hyperboloid
in Minkowski space:

\proclaim\@thm.CH1:Cor 5.1@.
The Christoffel dual of (part of) a tri-axial $1$-sheeted
hyperboloid (\:param-1hyp:) in Minkowski space,
$$
  u+jv=z\mapsto\alpha x(z)
  = {1\over\sn_qu\dn_qv}\pmatrix{
    a\cn_qu \cr -b\sn_qv\dn_qu \cr c\cn_qv \cr}
  \enspace{\sl with}\enspace
  q := \sqrt{b^2-c^2\over a^2-c^2},
$$
is the affine image of a timelike minimal surface,
more precisely,
up to homothety it is given by
$$
  z\mapsto\alpha x^\ast(z) := \Re\,\pmatrix{\hfill
   -{2a\over p}\arctan{\cn_p\over p\sn_p}\cr
   -j\,{2b\over pq}\arctan{q\cn_p\over p}\cr\hfill
   {2c\over q}\artanh q\sn_p\cr}(z)
  \enspace{\sl with}\enspace
  p := \sqrt{a^2-b^2\over a^2-c^2},
$$
where $z=u+jv\in\L$ is a para-complex variable, and
$\cn_p$ and $\sn_p$ denote the $\L$-Jacobi functions
obtained as $\L$-analytic extensions of the respective
real Jacobi elliptic functions.

To investigate the singularities of the curvature line net,
that is, umbilics of the $1$-sheeted hyperboloid,
we seek again the branch values of the differential equation
(\:orthoH1M:) resp (\:hypihyp:):
$$
  0 = p^2y^4 + 2(1+q^2)y^2 + p^2
  \enspace\Leftrightarrow\enspace
  y^2 = \cases{
    -({1\pm q\over p})^2
    = -({\sqrt{a^2-c^2}\pm\sqrt{b^2-c^2}
      \over\sqrt{a^2-b^2}})^2, \cr
    -({1\pm jq\over p})^2
    = -({\sqrt{a^2-c^2}\pm j\sqrt{b^2-c^2}
      \over\sqrt{a^2-b^2}})^2. \cr}
$$
Thus in order for the branch values to exist
(as para-complex numbers) we must have
$$
  {1\pm q\over p}
  = {\sqrt{a^2-c^2}\pm\sqrt{b^2-c^2}\over\sqrt{a^2-b^2}}
  \in i\R
  \enspace\Leftrightarrow\enspace\cases{
   c > a > b & or \cr
   b > a > c. & \cr}
$$
As the para-complex equation $y^2=1$ has four solutions,
$y=\pm 1$ and $y=\pm j$,
we obtain $16$ branch values in this case:
$$
  \pm {q\pm 1\over ip}, \enspace
  \pm {q\pm j\over ip}, \enspace
  \pm j{q\pm 1\over ip}, \enspace
  \pm j{q\pm j\over ip},
$$
four of which do not project to the hyperboloid since
$\|{q+j\over ip}\|^2={q+j\over ip}{q-j\over ip}=1$,
see \:fig.h1u:.
At the same time, the branch values yield the corners
of the range of the solution $y$ of (\:hypisol:),
as well as of the solutions $-y$ and $\pm jy$,
obtained from $y$ by symmetry.

 \pdfximage width 120pt {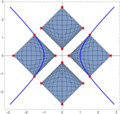}
  \xdef\hypmXd{\the\pdflastximage}
 \pdfximage width 3.5truein {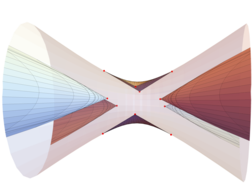}
  \xdef\hypmXs{\the\pdflastximage}
 \vglue -30pt\hbox to \hsize{\hfil
  \vtop{\vglue 56pt\pdfrefximage\hypmXd}%
  \hglue 32pt
  \vtop{\pdfrefximage\hypmXs}%
 \hfil}\vglue -0pt\centerline{%
  {\bf\@fig.h1u:Fig 6@.} Parametrizations of patches of
  a $1$-sheeted hyperboloid in $\R^{2,1}$ with umbilics}\par

This sharply contrasts the Euclidean case,
where the orthogonality condition leads to the differential
equation
$$
  y'^2 = \varrho\,\{ a^2(1+y^2)^2 + b^2(1-y^2)^2 + 4c^2y^2 \},
\eqno orthoH1E$$
cf (\:orthoEM:) and (\:orthoH2E:).
As a para-complex differential equations this does not admit
any branch values --- not surprisingly, as a $1$-sheeted
hyperboloid in Euclidean space cannot have any umbilics.

Following our earlier strategy we obtain an explicit
solution by replacing $a$ by $ia$ in (\:hypihyp:);
the corresponding curvature line parametrization of the
hyperboloid and its dual is then simply read off \:thm.CH1::

\proclaim\@thm.CH1E:Cor 5.2@.
A tri-axial $1$-sheeted hyperboloid
$
  1+({x_1\over a})^2=({x_2\over b})^2+({x_3\over c})^2
$
in Euclidean space, where without loss of generality $b<c$,
and its Christoffel dual admit curvature line parametrizations
by $\L$-analytic extensions of Jacobi elliptic functions:
$$
  \alpha x(z)
  = {1\over\sn_qu\dn_qv}\pmatrix{
    a\cn_qu \cr -b\sn_qv\dn_qu \cr c\cn_qv \cr}
  \enspace{\sl and}\enspace
  \alpha x^\ast(z) = \Re\,\pmatrix{\hfill
   -{2a\over p}\arctan{\cn_p\over p\sn_p}\cr
   -j\,{2b\over pq}\arctan{q\cn_p\over p}\cr\hfill
   {2c\over q}\artanh q\sn_p\cr}(z).
$$
where $z=u+jv\in\L$,
and $p=\sqrt{a^2+b^2\over a^2+c^2}$
and $q=\sqrt{-{b^2-c^2\over a^2+c^2}}$,
cf \:fig.h1e:.

\section Appendix. Reinbek's formulas

It has been known for long that quadrics are isothermic
surfaces, in fact, it seems that investigations on quadrics
were one of the main motivations to investigate isothermic
surfaces, cf [\:ca72:].
After obtaining the Christoffel duality for isothermic
surfaces in [\:ch67:] an obvious question was to determine
the Christoffel transforms of quadrics:
this was the topic of the PhD thesis [\:re86:].

As this thesis may not be very easy to obtain or read
we briefly list the results that are most relevant
to this paper:
using the well known conversion formulas
$$
  \artanh t = {1\over 2}\ln{1+t\over 1-t}
   \enspace{\rm and}\enspace
  \arctan t = {1\over 2i}\ln{1+it\over 1-it}
$$
we rewrite the representations that Reinbek obtained
for the Christoffel duals of quadrics in terms of the
usual elliptic coordinates in [\:re86:, \S\S6,7,8]:

{\parindent=2em

\item{(i)} for an ellipsoid 
  $\{ (x_1,x_2,x_3)\,|\,
   ({x_1\over a})^2+({y\over b})^2+({z\over c})^2=1\}$,
  with $a^2\geq t_1\geq b^2\geq t_2\geq c^2$,
  $$
    x(t_1,t_2) = \pmatrix{
     a\,\sqrt{(a^2-t_1)(a^2-t_2)\over(a^2-b^2)(a^2-c^2)} \cr
     b\,\sqrt{(b^2-t_1)(b^2-t_2)\over(b^2-a^2)(b^2-c^2)} \cr
     c\,\sqrt{(c^2-t_1)(c^2-t_2)\over(c^2-a^2)(c^2-b^2)} \cr
    }\enspace{\rm and}\enspace
    x^\ast(t_1,t_2) = \pmatrix{
     \sqrt{a^2\over(a^2-b^2)(a^2-c^2)}
      \artanh\sqrt{a^2-t_1\over a^2-t_2} \cr
     \sqrt{b^2\over(a^2-b^2)(b^2-c^2)}
      \arctan\sqrt{t_1-b\over b^2-t_2} \cr
     \sqrt{c^2\over(a^2-c^2)(b^2-c^2)}
      \artanh\sqrt{t_2-c^2\over t_1-c^2} \cr
    };
  $$

\item{(ii)} for a hyperboloid
  $\{ (x_1,x_2,x_3)\,|\,
   ({x_1\over a})^2+({y\over b})^2-({z\over c})^2=1\}$,
  with $a^2\geq t_1\geq b^2 > -c^2\geq t_2$,
  $$
    x(t_1,t_2) = \pmatrix{
     a\,\sqrt{(a^2-t_1)(a^2-t_2)\over(a^2-b^2)(a^2+c^2)} \cr
     b\,\sqrt{(b^2-t_1)(b^2-t_2)\over(b^2-a^2)(b^2+c^2)} \cr
     c\,\sqrt{(c^2+t_1)(c^2+t_2)\over(c^2+a^2)(c^2+b^2)} \cr
    }\enspace{\rm and}\enspace
    x^\ast(t_1,t_2) = \pmatrix{
     \sqrt{a^2\over(a^2-b^2)(a^2+c^2)}
      \artanh\sqrt{a^2-t_1\over a^2-t_2} \cr
     \sqrt{b^2\over(a^2-b^2)(b^2+c^2)}
      \arctan\sqrt{t_1-b^2\over b^2-t_2} \cr
     \sqrt{c^2\over(a^2+c^2)(b^2+c^2)}
      \arctan\sqrt{-{t_1+c^2\over c^2+t_2}} \cr
    };
  $$

\item{(iii)} for a hyperboloid
  $\{ (x_1,x_2,x_3)\,|\,
   ({x_1\over a})^2-({y\over b})^2-({z\over c})^2=1\}$,
  with $a^2>-b^2\geq t_1\geq -c^2\geq t_2$,
  $$
    x(t_1,t_2) = \pmatrix{
     a\,\sqrt{(a^2-t_1)(a^2-t_2)\over(a^2+b^2)(a^2+c^2)} \cr
     b\,\sqrt{(b^2+t_1)(b^2+t_2)\over(b^2+a^2)(b^2-c^2)} \cr
     c\,\sqrt{(c^2+t_1)(c^2+t_2)\over(c^2+a^2)(c^2-b^2)} \cr
    }\enspace{\rm and}\enspace
    x^\ast(t_1,t_2) = \pmatrix{
     \sqrt{a^2\over(a^2+b^2)(a^2+c^2)}
      \artanh\sqrt{a^2-t_1\over a^2-t_2} \cr
     \sqrt{b^2\over(a^2+b^2)(-b^2+c^2)}
      \artanh\sqrt{b^2+t_1\over b^2+t_2} \cr
     \sqrt{c^2\over(a^2+c^2)(-b^2+c^2)}
      \arctan\sqrt{-{t_1+c^2\over c^2+t_2}} \cr
    }.
  $$

\par}

\section Appendix. Lorentz analytic functions

Our analysis of the Christoffel duals of quadrics and the
associated zero mean curvature surfaces hinges crucially
on elliptic functions, as solutions of the occurring
elliptic differential equations.
To obtain the para-holomorphic analogues of the Jacobi
ellitpic functions we seek a suitable extension theorem
for para-holomorphic functions.

Following [\:tkl14:, Def 4.11(4)]
we say that a function $y:\L\supset U\to\L$ is
{\em $\L$-analytic\/} (or {\em Lorentz-analytic\/})
if it admits a power series expansion around every point
$z_0=u_0+jv_0\in U$,
$$
  y(u+jv) = \sum_{n\in\N}a_n((u-u_0)+j(v-v_0))^n.
$$
Any $\L$-analytic function is para-holomorphic in its domain,
in fact, has para-complex derivatives of any order,
which are given by term-by-term differentiation
of its power series [\:tkl14:, Thm 4.12].

Note that,
in contrast to the case of complex holomorphic functions,
para-holomorphicity does neither imply $\L$-analyticity
nor the existence of higher order derivatives,
cf [\:tkl14:, Expl 4.13].

In fact, any differentiable function $y:\R\supset I\to\R$,
regardless of higher order differentiability of $y$,
has a para-holomorphic extension
$$
  \tilde y(u+jv) :
  = y(u+v)\,{1+j\over 2} + y(u-v)\,{1-j\over 2}
$$
as
$
  \Re\tilde y(u,v) = {1\over 2}(y(u+v)+y(u-v))
$ and $
  \Im\tilde y(u,v) = {1\over 2}(y(u+v)-y(u-v))
$
clearly satisfy the para-complex version of the Cauchy-Riemann
equations,
$$
  (\Re\tilde y)_u = (\Im\tilde y)_v
   \enspace{\rm and}\enspace
  (\Im\tilde y)_u = (\Re\tilde y)_v.
$$
Note that this extension is not unique:
for example, $j\tilde y$ is para-holomorphic
with $j\tilde y|_I\equiv 0$,
hence any superposition of $\tilde y$ and $j\tilde y$
yields another extension of $y$;
in fact, any real differentiable function $y$ provides
a para-holomorphic extension $j\tilde y$ of $0$.

To see that the above construction yields a unique extension
of a real analytic function to an $\L$-analytic function,
first observe that
$$
  ({1\pm j\over 2})^2 = {1\pm j\over 2}
   \enspace{\rm and}\enspace
  ({1+j\over 2})({1-j\over 2}) = 0.
$$
Thus writing $u+jv=(u+v){1+j\over 2}+(u-v){1-j\over v}$
we learn that the extension of a real power series is
a para-complex power series since
$$
  \sum_{n\in\N}a_n(u+jv)^n
  = \sum_{n\in\N}a_n(u+v)^n{1+j\over 2}
  + \sum_{n\in\N}a_n(u-v)^n{1-j\over 2}.
$$
Uniqueness follows by comparison of coefficients:
if $\tilde y(u+jv)=\sum_{n\in\N}a_n(u+jv)^n$ restricts
to a real analytic function $y$ then all coefficients $a_n$
must be real, hence are uniquely determined by $y$.

Thus we have proved the following,
cf [\:tkl14:, Expl 3.2 and Prop 3.1(4)]:

\proclaim\@thm.laext:Lemma A.1@.
Every real analytic function $y:\R\supset I\to\R$ has
a unique $\L$-analytic extension
$$
  \tilde y:\L\supset\{u+jv\,|\,u\pm v\in I\}\to\L, \enspace
  u+jv\mapsto\tilde y(u+jv) :
   = y(u+v)\,{1+j\over 2} + y(u-v)\,{1-j\over 2}.
$$

Note that the $\L$-analytic trigonometric functions
of [\:tkl14:, Expl 3.3] are different from those that
are obtained by extension of the real trigonometric
functions.

Clearly, sums and products of the $\L$-analytic extensions
of real analytic functions are the $\L$-analytic extensions
of the corresponding sums resp products of the real functions,
and as
$$
  (y(u+v){1+j\over 2}+y(u-v){1-j\over 2})
  ({1\over y(u+v)}{1+j\over 2}+{1\over y(u-v)}{1-j\over 2})
  = 1
$$
the same holds true for quotients, where they are defined.
In fact, these properties also follow from the compatibility
of our $\L$-analytic extension with the composition,
$$
  \widetilde{f\circ y} = \tilde f\circ\tilde y
   \enspace{\rm since}\enspace
  (\Re\tilde y\pm\Im\tilde y)(u+jv) = y(u\pm v).
$$
the composition of the extensions of real analytic
functions 

As the derivative of an $\L$-analytic function is obtained
by term-by-term differentiation of its power series expansion
it is clear that the derivative of the $\L$-analytic extension
of a real analytic function is the $\L$-analytic extension
of its derivative.
However, a similar argument as above can be given,
that does not depend on the power series expansion:
$$
  \tilde y'(u+jv)
  = {1\over 2}
    ({\partial\over\partial u}+j{\partial\over\partial v})\,
    \tilde y(u+jv)
  = {y'+jy'\over 2}(u+v){1+j\over 2}
  + {y'-jy'\over 2}(u+v){1-j\over 2}
  = \widetilde{y'}(u+jv).
$$

\section Appendix. Jacobi elliptic functions

For the reader's convenience we collect some basic facts
about the Jacbi elliptic functions of pole type $n$:
in terms of the Jacobi amplitude function $\am_p$
these may be defined by
$$
  \cn_pz = \cos\am_pz, \enspace
  \sn_pz = \sin\am_pz, \enspace
  \dn_pz = \sqrt{1-p^2\sn_p^2z}.
$$
The Pythagorean laws are then a direct consequence of these
definitions:
$$
  1 = \cn_p^2 + \sn_p^2 = \dn_p^2 + p^2\sn_p^2;
$$
and with $\am_p'=\dn_p$ the characterizing elliptic
differential equations are readily derived:
$$
  \cn_p'^2 = (1-\cn_p^2)(q^2+p^2\cn_p^2), \enspace
  \sn_p'^2 = (1-\sn_p^2)(1-p^2\sn_p^2), \enspace
  \dn_p'^2 = (\dn_p^2-1)(q^2-\dn_p^2).
$$
The Jacobi elliptic functions are analytic functions,
naturally defined for complex arguments:
in fact, they are holomorphic functions from a torus
to the Riemann sphere since they are doubly periodic,
with the period lattice spanned by $4K_p$ and $4iK_q$,
where $K_p$ denotes the complete ellitpic integral
of the first kind.

To resort to real arguments the conversion formulae
of imaginary arguments,
$$
  \cn_p(iz) = {1\over\cn_qz}, \enspace
  \sn_p(iz) = i\,{\sn_qz\over\cn_qz}, \enspace
  \dn_p(iz) = {\dn_qz\over\cn_qz},
$$
are combined with the argument sum formulas
$$\matrix{
  \cn_p(u+v) = {
   \cn_pu\cn_pv-\sn_pu\dn_pu\sn_pv\dn_pv
    \over 1-p^2\sn_p^2u\sn_p^2v}, \hfill\cr
  \sn_p(u+v) = {
   \sn_pu\cn_pv\dn_pv+\sn_pv\cn_pu\dn_pu
    \over 1-p^2\sn_p^2u\sn_p^2v}, \hfill\cr
  \dn_p(u+v) = {
   \dn_pu\dn_pv-p^2\sn_pu\cn_pu\sn_pv\cn_pv
    \over 1-p^2\sn_p^2u\sn_p^2v}. \hfill\cr
}$$
Hence, for example,
$$
  e^{i\am_p(u+iv)}
  = {\cn_pu\cn_qv+i\sn_pu\dn_qv\over 1+\dn_pu\sn_qv}.
$$

Finally, the conversion formulas for inverse and imaginary
moduli allow to unify the treatment of hyperboloids
of different half axes configurations:
$$\matrix{
  \cn_{1/p}z = \dn_p{z\over p}, \hfill&
  \cn_{ip}z = {\cn_{p'}\over\dn_{p'}}{z\over q'}, \hfill\cr
  \sn_{1/p}z = p\sn_p{z\over p}, \hfill&
  \sn_{ip}z = q'{\sn_{p'}\over\dn_{p'}}{z\over q'}, \hfill\cr
  \dn_{1/p}z = \cn_p{z\over p}; \hfill&
  \dn_{ip}z = {1\over\dn_{p'}}{z\over q'}, \hfill\cr
}\enspace{\rm where}\enspace\cases{
  p'={p\over\sqrt{1+p^2}}, \cr
  q'={1\over\sqrt{1+p^2}}. \cr
}$$
Note that these conversions preserve reality of the
Jacobi elliptic functions for real arguments.

Considering the Jacobi elliptic functions as real analytic
functions we may use \:thm.laext: to obtain $\L$-analytic
extensions that are defined on a square torus:

\proclaim\@def.ljacobi:Def \& Cor A.2@.
The Jacobi elliptic functions have unique $\L$-analytic
extensions to $\L$, that we will refer to as
{\em$\L$-elliptic\/} or {\em$\L$-Jacobi functions\/}.
These $\L$-elliptic functions satify the same elliptic
differential equations as their real counterparts and they
are doubly periodic, with fundamental domains
$$
  \{ u+jv\,|\, -2K_p\leq u\pm v\leq 2K_p \}.
$$

Existence of unique $\L$-analytic extensions and their
periodicity follows directly from \:thm.laext:;
as the above elliptic differential equations characterizing
the Jacobi elliptic functions are (real) analytic,
$$
  0 = f(y',y)
   \enspace{\rm with}\enspace
  f\in C^\omega(\R^2),
$$
the extensions
of $y'$ and $y$ satisfy the same elliptic differential
equations as the original real functions do.
For example, the real solution
$$
  y(u) = e^{i\am_piu} = {\cn_qu\over 1+\sn_qu}
   \enspace{\rm of}\enspace
  y'^2 = {1\over 4}\{ p^2(1+y^4)+2(1+q^2)\,y^2 \}
$$
extends to a para-complex solution of the same elliptic
differential equation,
$$
  y(u+jv) = {\cn_qu-j\sn_qv\dn_qu\over\cn_qv+\sn_qu\dn_qv}.
$$

Note that $-y$ and $\pm jy$ provide other solutions
of the same elliptic differential equation.

\section Appendix. Implicit representations

While verification of the implicit representations of the
occurring zero mean curvature surfaces given in the text is
fairly straightforward, we provide some hints and formulas
that may facilitate the endeavour in this appendix.
First recall that three cases are discussed:

\vskip 1ex
\centerline{\vbox{\halign{\vrule height 12pt width 0pt depth 3pt
 \enspace#\enspace\vrule &
 \enspace#\enspace\vrule &
 \enspace#\enspace \cr
 Case &\hfil elliptic ode \hfil&\hfil solution \hfil\cr
 \noalign{\hrule}
 Ellipsoid/Scherk towers (\:elliell:) &\hfill
  $y'^2 = {(1-y^2)^2 - q^2(1+y^2)^2\over 4}$ &
  $y(z)={1\over i}\,e^{i\am_pz}$, \hfill\cr\hfill
 $2$-sheeted hyperboloid (\:hypiell:) &\hfill
  $-y'^2 = {(1+y^2)^2 - q^2(1-y^2)^2\over 4}$ &
  $y(z)=e^{i\am_pz}$, \hfill\cr\hfill
 $1$-sheeted hyperboloid (\:hypihyp:) &\hfill
  $y'^2 = {(1+y^2)^2 - q^2(1-y^2)^2\over 4}$ &
  $y(z) = e^{i\am_piz}$. \hfill\cr
}}}

Here
$p=\sqrt{{a^2-b^2\over a^2-c^2}}$
and
$q=\sqrt{1-p^2}=\sqrt{{b^2-c^2\over a^2-c^2}}$
denote the elliptic modulus and co-modulus of the
elliptic function $y$,
obtained from the half axes $a,b,c>0$ of the ellipsoid
or hyperboloid, respectively.
Note that, depending on the order of the half axes,
three cases can occur:
$$\matrix{
  \hfill{\rm(i)} & a>b>c &{\rm or}& c>b>a
   &{\rm yields}& p\in(0,1) &{\rm and}& q\in(0,1); \cr
  \hfill{\rm(ii)} & c>a>b &{\rm or}& b>a>c
   &{\rm yields}& p\in i\R &{\rm and}& q>1; \cr
  \hfill{\rm(iii)} & b>c>a &{\rm or}& a>c>b
   &{\rm yields}& p>1 &{\rm and}& q\in i\R. \cr
}$$

Using the differential equation of the elliptic function
$y$ the integrals in the Weierstrass representations,
 (\:intell:), (\:int2hyp:) and (\:int1hyp:),
are readily verified using
$$\matrix{
  (\artanh{2y\over p\,(1+y^2)})'
  &=& {2y'p\,(1-y^2)\over (1-y^2)^2-q^2(1+y^2)^2}, &
  (\arctan{2y\over p\,(1-y^2)})'
  &=& {2y'p\,(1+y^2)\over (1+y^2)^2-q^2(1-y^2)^2}; \cr
  (\artanh{2qy\over p\,(1-y^2)})'
  &=& {2y'pq\,(1+y^2)\over (1-y^2)^2-q^2(1+y^2)^2}, &
  (\arctan{2qy\over p\,(1+y^2)})'
  &=& {2y'pq\,(1-y^2)\over (1+y^2)^2-q^2(1-y^2)^2}; \cr
  (\artanh{q\,(1\pm y^2)\over (1\mp y^2)})'
  &=& {\pm 4yy'q\,\over (1\mp y^2)^2-q^2(1\pm y^2)^2}. \cr
}$$

We now investigate the three cases,
the ellipsoid in Euclidean space and 
the space- resp timelike hyperboloids in Minkowski geometry,
in turn.

\underbar{Ellipsoid/Scherk tower (\:intell:)}.
We set
$$\matrix{
  z_1 &:=& \hfill\int{1-y^2\over y'^2}\,dy
  &=&\,{2\over p}\,\artanh{2y\over p\,(1+y^2)}, \hfill\cr
  z_2 &:=& i\,\int{1+y^2\over y'^2}\,dy
  &=& {2i\over pq}\artanh{2qy\over p\,(1-y^2)} \hfill
  &=& {2\over pq}\,\arctan{2iqy\over p\,(1-y^2)}, \cr
  z_3 &:=& \int{2y\over y'^2}\,dy
  &=&\,{2\over q}\,\artanh{q\,(1+y^2)\over(1-y^2)}. \cr
}$$
Then we use
$$
  \cosh(z+\bar z)
  = {1+|\tanh z|^2\over|1-\tanh^2z|}
   \enspace{\rm and}\enspace
  \cos(z+\bar z)
  = {1-|\tan z|^2\over|1+\tan^2 z|}
$$
to compute
$$\matrix{\hfill
  \tanh^2{pz_1\over 2} = {4y^2\over p^2(1+y^2)^2}, &\hfill
  \cosh\Re pz_1 = {p^2|1+y^2|^2+4|y|^2
   \over|(1-y^2)^2-q^2(1+y^2)^2|}; \cr\hfill
  \tan^2{pqz_2\over 2} = {-4q^2y^2\over p^2(1-y^2)^2}, &\hfill
  \cos\Re pqz_2 = {p^2|1-y^2|^2-4q^2|y|^2
   \over|(1-y^2)^2-q^2(1+y^2)^2|}; \cr\hfill
  \tanh^2{qz_3\over 2} = {q^2(1+y^2)^2\over(1-y^2)^2}, &\hfill
  \cosh\Re qz_3 = {|1-y^2|^2+q^2|1+y^2|^2
   \over|(1-y^2)^2-q^2(1+y^2)^2|}. \cr
}$$
Without loss of generality we assume case (i) $a>b>c$ here,
so that $p,q\in(0,1)$.
Thus we obtain the implicit representation (\:implell:)
of the Scherk tower $x^\ast$ of \:thm.ST::
$$
  q^2\cosh\Re pz_1 + \cos\Re pqz_2 - p^2\cosh\Re qz_3 = 0.
$$

\underbar{$2$-sheeted hyperboloid (\:int2hyp:)}.
Here we set
$$\matrix{\hfill
  z_1 &:=&\hfill \int{1+y^2\over y'^2}\,dy
  &=&\,-{2\over p}\,\arctan{2y\over p\,(1-y^2)}, \hfill\cr
  z_2 &:=& i\,\int{1-y^2\over y'^2}\,dy
  &=& -{2i\over pq}\arctan{2qy\over p\,(1+y^2)} \hfill
  &=&  {2\over pq}\,\artanh{2qy\over ip\,(1+y^2)}, \cr
  z_3 &:=& \int{2y\over y'^2}\,dy
  &=&\,\phantom{-}{2\over q}\,\artanh{q\,(1-y^2)\over(1+y^2)}
  &=& {2i\over q}\,\arctan{q\,(1-y^2)\over i\,(1+y^2)}. \cr
}$$
Employing the same method as before,
two cases,
 (i) $a>b>c$ and (iii) $b>c>a$,
need to be considered,
as the $x_3$-axis is geometrically distinguished:
exchanging the $x_1$- and $x_2$-axes will interchange
the equations of cases (i) and (ii).
Thus we compute
$$\matrix{\hfill
  \tan^2{pz_1\over 2} =\phantom{-}
   {4y^2\over p^2(1-y^2)^2}, &\hfill
  \cos\Re pz_1 = {p^2|1-y^2|^2-4|y|^2
   \over|(1+y^2)^2-q^2(1-y^2)^2|}\rlap{;} \cr\hfill
  \tanh^2{pqz_2\over 2} =-{4q^2y^2\over p^2(1+y^2)^2}, &\hfill
  \cosh\Re pqz_2 = {p^2|1+y^2|^2+4q^2|y|^2
   \over|(1+y^2)^2-q^2(1-y^2)^2|}
  &{\rm if}& p<1, \cr\hfill
  \tan^2{pqz_2\over 2i} =\phantom{-}
   {4q^2y^2\over p^2(1+y^2)^2}, &\hfill
  \cos\Re{pqz_2\over i} = {p^2|1+y^2|^2+4q^2|y|^2
   \over|(1+y^2)^2-q^2(1-y^2)^2|}
  &{\rm if}& p>1; \cr\hfill
  \tanh^2{qz_3\over 2} =\phantom{-}
   {q^2(1-y^2)^2\over (1+y^2)^2}, &\hfill
  \cosh\Re qz_3 = {|1+y^2|^2+q^2|1-y^2|^2
   \over|(1+y^2)^2-q^2(1-y^2)^2|}
  &{\rm if}& p<1, \cr\hfill
  \tan^2{qz_3\over 2i} =-{q^2(1-y^2)^2\over (1+y^2)^2}, &\hfill
  \cos\Re{qz_3\over i} = {|1+y^2|^2+q^2|1-y^2|^2
   \over|(1+y^2)^2-q^2(1-y^2)^2|}
  &{\rm if}& p>1; \cr
}$$
and hence obtain implicit representations (\:impl2hyp:)
for the maximal surfaces of \:thm.CH2::
$$\matrix{
  q^2\cos\Re pz_1
  &+& \cosh\Re pqz_2
  &-& p^2\cosh\Re qz_3
  &=& 0
  &{\rm if}& p<1 \enspace{\rm and}\enspace q\in(0,1); \hfill\cr
  q^2\cos\Re pz_1
  &+& \cos\Re{pqz_2\over i}
  &-& p^2\cos\Re{qz_3\over i}
  &=& 0
  &{\rm if}& p>1 \enspace{\rm and}\enspace q\in i\R. \hfill\cr
}$$

\underbar{$1$-sheeted hyperboloid (\:int1hyp:)}.
Here we use para-complex coordinate functions
$$\matrix{\hfill
  z_1 &:=&\hfill \int{1+y^2\over y'^2}\,dy
  &=&\,{2\over p}\,\arctan{2y\over p\,(1-y^2)} &
  &=&\,-{2i\over p}\,\artanh{2iy\over p\,(1-y^2)}, \hfill\cr
  z_2 &:=& j\int{1-y^2\over y'^2}\,dy
  &=& {2j\over pq}\arctan{2qy\over p\,(1+y^2)} &
  &=& -{2ij\over pq}\artanh{2iqy\over p\,(1+y^2)}, \hfill\cr
  z_3 &:=& \int{2y\over y'^2}\,dy
  &=&\! -{2\over q}\,\artanh{q\,(1-y^2)\over(1+y^2)}. \cr
}$$
Note that we are now using the $\L$-analytic extensions
of real ellipsic functions, hence may need to replace
$\artanh$ by ${\rm arcoth}$ if necessary
---
however, this has no effect as we again use
$$
  \cosh(z+\bar z) = {1+\|\tanh z\|^2\over\|1-\tanh^2z\|}
   \enspace{\rm and}\enspace
  \cos(z+\bar z) = {1-\|\tan z\|^2\over\|1+\tan^2 z\|},
$$
together with $\tan jz=j\tan z$, $\tanh jz=j\tanh z$ and
$\|jy\|^2=-\|y\|^2$,
to compute
$$\matrix{\hfill
  \tan{pz_1\over 2} =\phantom{-} {2y\over p\,(1-y^2)}, &\hfill
  \cos\Re pz_1 =\phantom{-} {p^2\|1-y^2\|^2-4\|y\|^2
   \over\|(1+y^2)^2-q^2(1-y^2)^2\|}
  &{\rm if}& p\in\R, \hfill\cr\hfill
  \tanh{pz_1\over 2i} = -{2iy\over p\,(1-y^2)}, &\hfill
  \cosh\Re{pz_1\over i} = -{p^2\|1-y^2\|^2-4\|y\|^2
   \over\|(1+y^2)^2-q^2(1-y^2)^2\|}
  &{\rm if}& p\in i\R; \cr\hfill
  \tan{pqz_2\over 2} =\phantom{-} {2jqy\over p\,(1+y^2)}, &\hfill
  \cos\Re pqz_2 =\phantom{-} {p^2\|1+y^2\|^2+4q^2\|y\|^2
   \over\|(1+y^2)^2-q^2(1-y^2)^2\|}
  &{\rm if}& p\in\R, \hfill\cr\hfill
  \tanh{pqz_2\over 2i} = -{2ijqy\over p\,(1+y^2)}, &\hfill
  \cosh\Re{pqz_2\over i} = -{p^2\|1+y^2\|^2+4q^2\|y\|^2
   \over\|(1+y^2)^2-q^2(1-y^2)^2\|}
  &{\rm if}& p\in i\R; \cr\hfill
  \tanh{qz_3\over 2} = -{q\,(1-y^2)\over(1+y^2)}, &\hfill
  \cosh\Re qz_3 =\phantom{-} {\|1+y^2\|^2+q^2\|1-y^2\|^2
   \over\|(1+y^2)^2-q^2(1-y^2)^2\|}\rlap{.} \cr
}$$
Thus we obtain implicit representations (\:impl1hyp:)
of the timelike minimal surfaces of \:thm.CH1::
$$\matrix{
  q^2\cos\Re pz_1
  &+& \cos\Re pqz_2
  &-& p^2\cosh\Re qz_3
  &=& 0
  &{\rm if}& p\in(0,1) \enspace{\rm and}\enspace q\in(0,1); \hfill\cr
  q^2\cosh\Re{pz_1\over i}
  &+& \cosh\Re{pqz_2\over i}
  &+& p^2\cosh\Re qz_3
  &=& 0
  &{\rm if}& p\in i\R \enspace{\rm and}\enspace q\in(1,\infty).
   \hfill\cr
}$$
Note that the third case, (iii) $b>c>a$, can be neglected
since the $x_1$-axis is the timelike axis so that this case
can be obtained from case (i) by exchanging the $x_2$- and
$x_3$-axes.

\section References

\message{References}
\bgroup\frenchspacing\parindent=2em

\refitem ca72
 A Cayley: 
 {\it On the surfaces divisible into squares by their curves
  of curvature\/};
 Proc London Math Soc 4, 8--9 \& 120--121 (1872)

\refitem ch67
 E Christoffel:
 {\it Ueber einige allgemeine Eigenschaften der
  Minimumsfl\"achen\/};
 Crelle's J 67, 218--228 (1867)

\refitem co49
 J Cockle:
 {\it On a new Imaginary in Algebra\/};
 Phil Mag 34, 37--47 (1849)

\refitem tkl14
 L di Terlizzi, J Konderak, I Lacirasella:
 {\it On differentiable functions over Lorentz numbers and
  their geometric applications\/};
 Differ Geom Dynam Syst 16, 113--139 (2014)

\refitem fkk12
 S Fujimori, Y Kim, S-E Koh, W Rossman, H Shin,
  H Takahashi, M Umehara, K Yamada, S-D Yang:
 {\it Zero mean curvature surfaces in ${\bb L}^3$
  containing a light-like line\/};
 C R 350, 975--978 (2012)

\refitem fkk15
 S Fujimori, Y Kim, S-E Koh, W Rossman, H Shin,
  M Umehara, K Yamada, S-D Yang:
 {\it Zero mean curvature surfaces in Lorentz-Minkowski
  $3$-space and $2$-dimensional fluid mechanics\/};
 Math J Okayama Univ 57, 173--200 (2015)

\refitem haha04
 A Harkin, J Harkin:
 {\it Geometry of generalized complex numbers\/};
 Math Mag 77, 118--129 (2004)

\refitem imdg
 U Hertrich-Jeromin:
 {\it Introduction to M\"obius differential geometry\/};
 London Math Soc Lect Note Series 300,
  Cambridge Univ Press, Cambridge (2003)

\refitem ka89
 H Karcher:
 {\it Construction of minimal surfaces\/};
 Surveys in Geometry 1--96, Tokyo Univ (1989)

\refitem ko83
 O Kobayashi:
 {\it Maximal surfaces in the $3$-dimensional Minkowski space
  $L^3$\/};
 Tokyo J Math 6, 297--309 (1983)

\refitem ko99
 J Konderak:
 {\it A Weierstrass representation theorem for
  Lorentz surfaces\/};
 Complex Variables, Theory Appl 50, 319--332 (2005).

\refitem ma04
 M Magid:
 {\it Lorentzian isothermic surfaces and Bonnet pairs\/};
 Ann Polon Math 83, 129--139 (2004)

\refitem ma05
 M Magid:
 {\it Lorentzian isothermic surfaces in $R^n_j$\/};
 Rocky Mountain J Math 35, 627--640 (2005)

\refitem ni89
 J Nitsche:
 {\it Lectures on minimal surfaces\/};
 Cambridge Univ Press, Cambridge (1989)

\refitem re86
 K Reinbek:
 {\it Ueber diejenigen Fl\"achen, auf welche die Fl\"achen
  zweiten Grades durch parallele Normalen conform abgebildet
  werden\/};
 PhD thesis, Georg-Augusts-Universit\"at G\"ottingen (1886)

\refitem sc35
 H Scherk:
 {\it Bemerkungen \"uber die kleinste Fl\"ache innerhalb
  gegebener Grenzen\/};
 Crelle's J 13, 185--208 (1835)

\refitem sm04
 B Smyth:
 {\it Soliton surfaces in the mechanical equilibrium
  of closed membranes\/};
 Commun Math Phys 250, 81--94 (2004)

\egroup
\vskip3em
\bgroup\fn[cmr7]\baselineskip=8pt
\def\addwd{\hsize=.36\hsize}
\def\udo{\vtop{\addwd
 U Hertrich-Jeromin\\
 Vienna University of Technology\\
 Wiedner Hauptstra\ss{}e 8--10/104\\
 A-1040 Vienna (Austria)\\
 Email: udo.hertrich-jeromin@tuwien.ac.at
 }}
\def\fujimori{\vtop{\addwd
 S Fujimori\\
 Department of Mathematics\\
 Okayama University\\
 Okayama 700-8530 (Japan)\\
 Email: fujimori@math.okayama-u.ac.jp
 }}
\def\kokubu{\vtop{\addwd
 M Kokubu\\
 Department of Mathematics\\
 Tokyo Denki University\\
 Tokyo 120-8551 (Japan)\\
 Email: kokubu@cck.dendai.ac.jp
 }}
\def\umehara{\vtop{\addwd
 M Umehara\\
 Department of Mathematical and Computing Sciences\\
 Tokyo Inistitute of Technology\\
 Tokyo 152-8552 (Japan)\\
 Email: umehara@is.titech.ac.jp
 }}
\def\yamada{\vtop{\addwd
 K Yamada\\
 Department of Mathematics\\
 Tokyo Inistitute of Technology\\
 Tokyo 152-8551 (Japan)\\
 Email: kotaro@math.titech.ac.jp
 }}
\hbox to \hsize{\hfil \fujimori \hfil \udo \hfil}\vskip 3ex
\hbox to \hsize{\hfil \kokubu \hfil \umehara \hfil}\vskip 3ex
\hbox to \hsize{\hfil \yamada \hfil}
\egroup
\end{document}                    